\documentclass[11pt]{article}
\usepackage{amsfonts,enumerate,hhline,amsthm}
\def\TheVersion{1.0}
\def\TheDate{November 24, 2003}
\def\TheTitle{3-transposition Groups of Symplectic Type\\
and\\
Vertex Operator Algebras\\
}
\newtheorem{lemma}{Lemma}[subsection]
\newtheorem{prop}[lemma]{Proposition}

\newtheorem{remark}[lemma]{\sl Remark}
\newtheorem{note}[lemma]{\sl Note}

\newtheorem{theorem}[lemma]{Theorem}
\newtheorem{corollary}[lemma]{Corollary}
\newcommand{\Prop}{\begin{prop}}
\newcommand{\EndProp}{\end{prop}}
\newcommand{\Remark}{\begin{remark}}
\newcommand{\EndRemark}{\end{remark}}
\newcommand{\Note}{\begin{note}}
\newcommand{\EndNote}{\end{note}}
\newcommand{\Lemma}{\begin{lemma}}
\newcommand{\EndLemma}{\end{lemma}}
\newcommand{\Theorem}{\begin{thm}}
\newcommand{\EndTheorem}{\end{thm}}

\catcode`\@=11

\@addtoreset{equation}{section}
\makeatother
\catcode`\@=\active
\def\EQ#1{{\def\normalbaselines{\baselineskip15pt \lineskip6pt \lineskiplimit4pt}\,\vcenter{\normalbaselines\mathsurround=0pt\ialign{\hfil$\displaystyle##$\hfil&&$\,$$\displaystyle##$\hfil\crcr\mathstrut\crcr\noalign{\kern-\baselineskip}#1\crcr\mathstrut\crcr\noalign{\kern-\baselineskip}}}\,}}
\def\tildeb#1{{\bf\tilde{\fam=-1 #1}}}%
\def\widetildeb#1{{\bf\widetilde{\fam=-1 #1}}}%
\def\barb#1{{\bf\bar{\fam=-1 #1}}}%
\newcommand{\FF}{{\mathbb F}}
\newcommand{\ZZ}{{\mathbb Z}}
\newcommand{\RR}{{\mathbb R}}

\newcommand{\CC}{{\mathbb C}}

\newcommand{\Lc}{{\cal L}}
\newcommand{\id}{{\rm id}}
\newcommand{\1}{{\bf 1}}
\newcommand{\ccc}{\gamma}
\newcommand{\cw}{\delta}
\newcommand{\Aut}{\mathop{\rm Aut}}
\newcommand{\Sym}{\mathop{\rm S}\nolimits}
\newcommand{\Symm}{\mathop{\rm Sym}\nolimits}
\newcommand{\Orth}{\mathop{\rm O}\nolimits}
\newcommand{\Sp}{\mathop{\rm Sp}\nolimits}
\newcommand{\compo}{{\mathchoice{\scriptstyle\circ}{\scriptstyle\circ}{\scriptscriptstyle\circ}{\scriptscriptstyle\circ}}}
\newcommand{\split}{{\,:\,}}
\newcommand{\VA}[1]{{\langle#1\rangle^{}_{\rm VOA}}}
\newcommand{\EE}{F}
\newcommand{\GL}{\rm GL}
\newcommand{\PSL}{\rm PSL}
\newcommand{\Col}{{\rm rel\,}}
\newcommand{\Point}{{(1)}}
\newcommand{\Line}{{(2)}}

\begin{document}
\title{
\vskip-8ex{\small\hfill Ver.\ \TheVersion,\ \ \TheDate
\kern-9ex \\}\vskip6ex
\Large
\TheTitle}
\author{
\cr
Atsushi Matsuo
\cr\cr
{\normalsize\rm Graduate School of Mathematical Sciences}\cr
{\normalsize\rm University of Tokyo, Komaba, Tokyo 153-8914, Japan}
}
\date{}
\maketitle
{\abstract\noindent
The 3-transposition groups that act on a vertex operator algebra in the way described by Miyamoto in \cite{Mi1} are classified under the assumption that the group is centerfree and the VOA carries a positive-definite invariant Hermitian form.
This generalizes and refines the result of Kitazume and Miyamoto \cite{KM}.
Application to a similar but different situation is also considered in part by a slight generalization of the argument.
}
\section*{Introduction}
A (finite) $3$-transposition group is a pair $(G,D)$ of a finite group $G$ and a normal set $D$ of involutions in $G$ such that $G$ is generated by $D$ and if $a,b\in D$ then the order of $ab$ is either $1$, $2$ or $3$.
We sometimes say that a group $G$ is a 3-transposition group without mentioning the set $D$ as far as the set $D$ is easily recognized or uniquely determined. In this paper, the group $\Sym_6$, the symmetric group of degree $6$, and its extensions are the only cases where the set $D$ is not uniquely determined.
A subgroup $H$ of $G$ is called a $D$-subgroup if the pair $(H,D\cap H)$ is a 
3-transposition group.

$3$-transposition groups are first considered and studied by Fischer \cite{Fi}, who discovered the sporadic Fischer groups on the way of classification of $3$-transposition groups under the condition that the group is centerfree and reduced.
The classificaton problem of $3$-transposition groups has been further  investigated by a number of researchers; See the paper \cite{CH} and references therein.
Those investigations seem to be aiming at removing the assumptions of finiteness and reducedness of the groups under classification by improving or replacing some group-theoretical arguments.

The study of centerfree 3-transposition groups is actually equivalent to that of a certain combinatorial structure on the set $D$ called the structure of a Fischer space. Roughly speaking, a Fischer space is a partial linear space which is made by glueing the dual affine plane of order $2$ and the affine plane of order $3$.
When the affine plane of order $3$ does not occur, then the corresponding 3-transposition group is said to be of symplectic type.

\vskip2ex

On the other hand, looking at some general features of the `Griess algebra' of a vertex operator algebra (VOA) concerned with actions the Virasoro algebra of the central charge $1/2$, Miyamoto \cite{Mi1} found that $3$-transposition groups arise as a group of automorphisms of a VOA under certain circumstance.

To be more precise, let $V=V_0\oplus V_1\oplus\cdots$ be a VOA such that $\dim V_0=1$ and $V_1=0$.
Let $e$ be a vector in the `Griess algebra' $B=V_2$ such that it generates the Virasoro VOA $L(1/2,0)$ of the central charge $1/2$, of which the irreducible modules are only $L(1/2,0)$, $L(1/2,1/2)$ and $L(1/2,1/16)$. 
If $V$ does not have components isomorphic to $L(1/2,1/16)$ then the irreducible decomposition of $V$ with respect to the $L(1/2,0)$ gives rise to an involution denoted by $\sigma_e$ which inverts the vectors belonging to components isomorphic to $L(1/2,1/2)$ leaving the components isomorphic to $L(1/2,0)$ uneffected.
Then for two such vectors $e$ and $f$ the order of $\sigma_e\sigma_f$ must be $1$, $2$ or $3$.

Since we are mainly interested in its action on the Griess algebra, we assume that the VOA $V$ is generated by such Virasoro vectors as above.
This assumption implies that the resulting 3-transposition group is centerfree.
Let us say that a centerfree $3$-transposition group $(G,D)$ is $(1/2,1/2)$-realizable by $V$ if it can be obtained in such a way as described above from the VOA $V$, where $D$ is the set of automorphisms of the form $\sigma_e$ and $G=\langle D\rangle$.

A particular case of such a VOA is obtained by constructing the VOA associated with a doubly even binary code, which was also considered by Miyamoto.
In \cite{KM}, Kitazume and Miyamoto gave a list of candidates of $3$-transposition groups which might be $(1/2,1/2)$-realizable by code VOAs and showed that the groups in the list are in fact $(1/2,1/2)$-realizable by  code VOAs except one.
They assumed that their code VOAs carry a positive-definite invariant Hermitian form.

\vskip2ex

The primary purpose of this paper is to give the classification of the
 $(1/2,1/2)$-realizable 3-transposition groups without assuming 
that the VOA is a code VOA.

\vskip2ex

Our main result is the following: A 3-transposition group is $(1/2,1/2)$-realizable by a VOA having a positive-definite invariant Hermitian form if and only if it is the direct product of a finite number of  $D$-subgroups of the groups $\EE^2\split \Sym_n$ and $\Orth^+_{10}(2)$, where $\EE=2^{2m}$ if $n=2m+1,2m+2$ (Theorem \ref{MainTheorem}).
The groups $\EE^2\split \Sym_n$ and $\Orth^+_{10}(2)$ are realized by the VOA $V^+_{\sqrt{2}D_n}$ and by $V^+_{\sqrt{2}E_8}$ respectively, their 3-transposition subgroups are realized by taking subVOAs, and the direct product is realized by taking the tensor product of VOAs.
Consequently, we see that the exceptional one in the list of Kitazume-Miyamoto is actually not $(1/2,1/2)$-realizable by a VOA.

In deriving the main result mentioned above, the key is the following  arguments, both of which use the existence of a positive-definite invariant Hermitian form on our VOA.
One is to consider the decomposition not only with respect to $L(1/2,0)$ but also to $L(7/10,0)$.
Then we see that the affine plane of order $3$ does not occur in the Fischer space associated with an $(1/2,1/2)$-realizable $(G,D)$ so that our $(G,D)$ is a 3-transposition group of symplectic type. 
We can now use the classification of such groups by J.I.~Hall \cite{Ha1},\cite{Ha2}.

The other is to consider the adjacency matrix of the collinearity graph of the Fischer space, called the diagram, associated with the set $D$.
Then we see that if $(G,D)$ is $(1/2,1/2)$-realizable\ by a VOA with a positive-definite invariant Hermitian form then the least eigenvalue of the diagram must be greater than or equal to $-8$.
Therefore, looking at the parameters of the $3$-transposition groups,   we eliminate cases from the list of $(1/2,1/2)$-realizable groups.

We then check that the cases that passed the above-mentioned tests are indeed realized by VOAs.
Each of those VOAs is actually obtained as a subVOA (with a possibly different conformal vector) of $V^+_{\sqrt{2}R}$, where $R$ is a root system.
The structure of the Griess algebra of the latter VOA was described by Dong et al.{} \cite{DLMN}.

\vskip2ex

We would like to mention that the arguments above can be partly generalized.
Namely we may replace the concept of $(1/2,1/2)$-realizable groups by that of $(\ccc,\cw)$-realizable groups if for the central charge $c=\ccc$ and the conformal weight $h=\cw$ the pair $L(\ccc,0)$ and $L(\ccc,\cw)$ have the nice fusion property as that of $L(1/2,0)$ and $L(1/2,1/2)$.

Therefore, we will start by investigating a general properties of the Griess algebra of a VOA concerning 3-transposition groups which are $(\ccc,\cw)$-realizable for general $\ccc,\cw$.
Then we will apply these arguments to the case $\ccc=1/2$, $\cw=1/2$ and will derive our main result.
We will give some statements for a few other $\ccc,\cw$ as well; However, as we do not know the general scheme of constructing such VOAs,  we cannot give the classification of those groups in the present paper.

\paragraph{Acknowledgement}
The author is grateful to A.A.~Ivanov for discussion, especially for his instruction on 3-transposition groups as well as on graph theory which was crucial in completing this work.
The author is also grateful to M.~Kitazume, J.~Duncan and C.-H.~Lam for discussion.
The author thanks S.-J.~Cheng, R.-B.~Zhang, C.-H.~Lam, N.~Lam, M.~Primc and D.~Adamovic for hospitality and conversation.
This work was partly supported by Grants-in-Aid for Scientific Research No.~14740005, Japan Society for the Promotion of Science.
\vskip7ex
\section{Griess Algebra of Vertex Operator Algebras}
We will describe some general features of vertex operator algebras, especially on the construction of their subalgebras.
\subsection{Vertex Operator Algebras}
\label{Vertex Operator Algebras}
We start by recalling the notion of a vertex operator algebra (VOA) just in order to fix our notation and terminology.
See e.g.{} \cite{Bor},\cite{FLM},\cite{FHL} and \cite{MN} for accounts on VOAs and their modules.

\vskip2ex

Consider a VOA $V$ over the field $\CC$ of complex numbers.
By definition, $V$ is a graded vector space over $\CC$ which is given a correspondence that assigns for each $u\in V$ a `field' $Y(u,z)$ of the form
\begin{equation}
Y(u,z)=\sum_{n\in\ZZ}u_{(n)}z^{-n-1},
\end{equation}
where the coefficients $u_{(n)}$ are operators acting on $V$.
We may view this assignment as giving the set of binary operations $(u,v)\mapsto u_{(n)}v$ on $V$ for all $n\in\ZZ$.
Besides this, a vector $\1$ callded the vacuum vector and a vector $\omega$ called the conformal vector (Virasoro element) are specified in $V$.
These structures satisfy number of properties, some of which are taken as the set of axioms for a VOA.

Among other things, the grading of the vector space $V$ is given by the eigenvalues of the operator $L_0=\omega_{(1)}$ and the graded piece $V_k$ of degree $k$ is supposed by definition to be finite-dimensional.
The vacuum vector $\1$ belongs to $V_0$ and the conformal vector $\omega$  to $V_2$.
The operator $L_0$ is a part of the set of operators $L_n=\omega_{(n+1)}$, $n\in\ZZ$.
The vector $\omega$ is supposed to satisfy the properties that the operators $L_n$ form an action of the Virasoro algebra and that the equality  $u_{(-2)}\1=\omega_{(0)}u$ holds for all $u\in V$.

The central charge of the Virasoro action generated by the conformal vector $\omega$ as above is called the {\it rank}\/ of the VOA $V$, which we will denote by the symbol $c^{}_V={\rm rank}\,V$.
It is determined by the relation $2\omega_{(3)}\omega=c^{}_V\1$.
Note that $V$ as a vector space is usually infinite-dimensional.

The operators $u_{(n)}$ for $u\in V$ and $n\in\ZZ$ satisfy a set of equalities  (\cite{Bor}), called the {\it Cauchy-Jacobi identity}\/ (\cite{FLM}) or the Borcherds identity:
\begin{equation}
\sum_{i=0}^\infty
{p\choose i}(u_{(r+i)}v)_{(p+q-i)}w
=\sum_{i=0}^\infty(-1)^i{r\choose i}(
u_{(p+r-i)}(v_{(q+i)}w)-(-1)^r v_{(q+r-i)}(u_{(p+i)}w))
\label{CJ}
\end{equation}
for all $p,q,r\in \ZZ$ and $u,v,w\in V$.

\vskip2ex

For instance, the irreducible highest weight module $L(c,0)$ over the Virasoro algebra with the central charge $c$ carries a unique structure of a VOA for which the vacuum vector $\1$ is the highest weight vector and the conformal vector is the vector $L_{-2}\1$.
\subsection{Invariant Hermitian form}
\label{Invariant Hermitian form}
An {\it invariant}\/ bilinear form on a VOA $V$ is a bilinear form $(\ |\ ):V\times V\rightarrow \CC$ such that
\begin{equation}
(Y(u,z)v|w)=(v|Y(e^{zL_1}(-z^{-2})^{L_0}u,z^{-1})w)
\label{inv}
\end{equation}
for any $u,v,w\in V$  (\cite{Bor}).
An invariant bilinear form is always symmetric (\cite{FHL}).
We normalize the form so that $(\1|\1)=1$.

\vskip2ex

Now consider the field $\RR$ of real numbers.
We may consider VOAs over $\RR$ as well as those over $\CC$.
A {\it real form}\/ of a VOA $V$ is a subspace $V_\RR\subset V$ over $\RR$ containing $\1$ and $\omega$ and satisfying $V_\RR\otimes_\RR\CC=V$ such that $V_\RR$ has a structure of a VOA over $\RR$ by restricting that of $V$.
Consider a pair of a real form $V_\RR$ and an invariant bilinear form $(\ |\ )$ on $V_\RR$ valued in $\RR$, which is extended to a Hermitian form $\langle\ |\ \rangle$ on $V$.
It satisfies
$\langle Y(u,z)v|w\rangle=\langle v|Y(e^{zL_1}(-z^{-2})^{L_0}a,z^{-1})w\rangle
$
for all $u\in V_\RR$, $v,w\in V$.
We will call such a {\it pair}\/ of $V_\RR$ and $\langle\ |\ \rangle$ an {\it invariant Hermitian} \/ form on the VOA $V$.

\vskip2ex

In the rest of the paper, we assume the following condition.
\begin{description}
\item{(I)}
$V$ has a positive-definite invariant Hermitian form.
\end{description}
Since $L_1\omega=0$, we have $\langle L_0v|w\rangle=\langle v|L_0w\rangle$ by the invariance (\ref{inv}).
Hence $\langle V_m|V_n\rangle=0$ if $m\ne n$.
Therefore for any graded subspace $U$ we may consider the orthogonal complement $U^\perp$ in $V$, which satisfies $U^\perp\cap V_n=(U\cap V_n)^\perp$ for all $n\geq 0$.

For any subspace $A$ of $V$, we will denote the intersection with the real form $V_\RR$ by
\begin{equation}
A_\RR=A\cap V_\RR.
\end{equation}
For example, $B_\RR$ means the corresponding real form of the Griess algebra $B$.

\vskip2ex

As $\langle\ |\ \rangle$ agrees with $(\ |\ )$ on the real form $V_\RR$, we will rather use $(\ |\ )$ when we will be dealing with vectors spanned by vectors belonging to $V_\RR$.
In fact, we could work over the field $\RR$ throughout the paper.
\subsection{Idempotents in the Griess algebra}
\label{Idempotents in the Griess algebra}
From now on, we only consider VOAs with the following property:
\begin{description}
\item{(II)}
$V=\bigoplus_{n=0}^\infty V_n$, $V_0=\CC\1$ and $V_1=0$.
\end{description}
Then the degree $2$ subspace $B=V_2$ has a structure of a commutative nonassociative algebra by setting $x\cdot y=x_{(1)}y$, for which the conformal vector $\omega$ is twice the unity.
The algebra $B$ has a unique symmetric bilinear form $(\ |\ )_B$ such that $(x|y)_B\1=x_{(3)}y$ which is invariant in the sense that $(a\cdot x|y)_B=(x|a\cdot y)_B$ for all $a,x,y\in B$.
These properties follow from the Cauchy-Jacobi identity (\ref{CJ}).

Let $(\ |\ )$ be an invariant bilinear form on $V$ normalized so that $(\1|\1)_V=1$.
The assumption (II) implies that $L_1V_2=0$.
Therefore, by the invariance (\ref{inv}), we have
\begin{equation}
(a_{(n)}v|w)=(v|a_{(2-n)}w)
\end{equation}
for all $a\in V_2$ and $v,w\in V$.
Therefore, for $x,y\in B$, we have $(x|y)_B=(x_{(3)}y|\1)=(y|x_{(-1)}\1)=(y|x)=(x|y)$.

Such an algebra $B$ equipped with the form $(\ |\ )$ is called the Griess algebra of $V$, for the algebra $B$ in the case of the moonshine module $V^\natural$ is a variant of the algebra originally constructed by Griess.

\vskip2ex

Let $V$ be a VOA satisfying (II) and consider the Griess algbera $B=V_2$.
We will say by abuse of terminology that a nonzero vector $e\in B=V_2$ satisfying $e\cdot e=2e$ an {\it idempotent}\/ of the Griess algebra $B$

Let $e$ be an idempotent of the Griess algbera.
Then it is obvious from the grading of $V$ that $e_{(n)}e=0$ if $n\geq 4$ or $n=2$.
This means that the field $T_e(z)=Y(e,z)$ satisfies the Virasoro operator product expansion
\begin{equation}
T_e(z)T_e(w)\sim \frac{c_e/2}{(z-w)^4}+\frac{2T_e(w)}{(z-w)^2}+\frac{\partial T_e(w)}{z-w},
\end{equation}
so that the modes $L^e_n=e_{(n+1)}$ of $T_e(z)$ are subject to the Virasoro commutation relation with the central charge $c_e=2(e|e)$:
\begin{equation}
[L^e_m,L^e_n]=(m-n)L^e_{m+n}+\frac{m^3-m}{12}\delta_{m+n,0}\,c_e\,\id.
\end{equation}
By abuse of terminology, we will call an eigenvalue of the operator $L^e_0$ acting on $B$ an {\it eigenvalue of the idempotent}\/ $e$.
\subsection{Vertex operator subalgebras}
Throughout the paper, we mean by a vertex operator subalgebra of a VOA $V$, or rather by a {\it subVOA}\/, a graded subspace $U$ which has a structure of a VOA such that the binary operations and the grading of $U$ agrees with the restriction of those of $V$ and that $U$ and $V$ share the same vacuum vector.
However, we do {\it not}\/ assume that they share the same conformal vector.
When they do, we will call $U$ a {\it full}\/ subVOA.

Let $S$ be any subset of $V$. Consider the subspace $\VA{S}$ generated by $S$, i.e., the span of the elements of $V$ of the form
\begin{equation}
a^1_{(n_1)}a^2_{(n_2)}\cdots a^k_{(n_k)}\1
\end{equation}
where $k$ is a nonnegative integer, $a^1,\ldots,a^k\in S$ and $n_1,\ldots,n_k\in \ZZ$.
Then the subspace $\VA{S}$ is in fact closed under the binary operations.
If it owns a conformal vector with appropriate properties then it is a subVOA in our sense.

\begin{lemma}
\sl
Let $V$ be a VOA and let $U$ be a graded subspace of $V$ containing the vacuum vector $\1$ of $V$ such that $U_2$ is nonzero.
Suppose that $U$ is closed under the binary operations of $V$ and that there exists a subspace $W$ such that $V=U\oplus W$ and $U_{(n)}W\subseteq W$.
Then  $U$ has a unique conformal vector that gives $U$ a structure of a subVOA of $V$.
\label{2.2}
\end{lemma}
\begin{proof}
Let $\omega=\xi+\eta$ be the decomposition of the conformal vector $\omega$ of $V$ with respect to $V=U\oplus W$.
We will show that the vector $\xi$ has the desired properties.
First note that $\xi,\eta\in V_2$, so 
$2\xi=\omega\cdot \xi=(\xi+\eta)\cdot \xi=\xi\cdot \xi+\xi\cdot \eta$.
Since $\xi\in U$, $\xi\cdot \xi\in U$ and $\xi\cdot \eta\in W$,  we have $\xi\cdot \xi=2\xi$ and $\xi\cdot \eta=0$.
Thus $\xi$ is an idempotent of $B\cap U$, which gives rise to an action of the Virasoro algebra on $U$.
Furthermore, for any $u\in U\cap V_n$, we have $nu=L_0u=\omega_{(1)}u=(\xi+\eta)_{(1)}u=\xi_{(1)}u+\eta_{(1)}u$.
Since $nu$ and $\xi_{(1)}u$ belong to $U$ and $\eta_{(1)}u$ belongs to $W$, we see that $L^{\xi}_0u=\xi_{(1)}u=nu$.
Hence the grading of $U$ with respect to $L^\xi_0$ agrees with the restriction of the grading of $V$.
Finally, for any $u\in U$ we have $u_{(-2)}\1=\omega_{(0)}u={\xi}_{(0)}u+{\eta}_{(0)}u$.
However, since $u_{(-2)}\1\in U$, ${\xi}_{(0)}u\in U$ and ${\eta}_{(0)}u\in W$, we see that $u_{(-2)}\1={\xi}_{(0)}u$.
The uniqueness is obvious.
\end{proof}

\begin{lemma}
\sl
Let $V$ be a VOA satisfying the condition {\rm (I)}.
Let $U$ be a graded subspace of $V$ closed under the binary operations such that $U=U_\RR\otimes_\RR\CC$.
Then the orthogonal complement $W$ of $U$ with respect to the Hermitian form satisfies $U_{(n)}W\subseteq W$ for all $n\in \ZZ$.
\end{lemma}
\begin{proof}
From the assumptions it follows that $U_\RR$ is closed under the binary operations. So we have $(U|u_{(n)}w)\subseteq (U|w)=0$ for any $u\in U_\RR$ and $w\in W$.
This shows that $U_{(n)}W\subseteq W$ because $U=U_\RR\otimes_\RR\CC$.
\end{proof}

\vskip2ex

The lemmas above immediately imply the following.

\begin{prop}
\sl
Let $V$ be a VOA satisfying the conditions {\rm (I)} and {\rm (II)}.
Then, for any nontrivial subset $S$ of $B_\RR$, the space $\VA{S}$ generated by $S$ has a unique structure of a subVOA of $V$.
\end{prop}

Now let us consider the case when the subVOA is generated by a subalgebra of the Griess algebra.
The next lemma gives a sufficient condition for the Griess algebra of the subVOA to be equal to the generating subalgebra.
\begin{lemma}
\sl
Let $V$ be a VOA satisfying {\rm (I)} and {\rm (II)} and let $A$ be a subalgebra of $B$.
If there exists a subset $T\subset B_\RR$ such that $A=\{a\in B\,|\,a\cdot T=0\}$, then $V_2\cap \VA{A}=A$.
\label{lemma: U2}
\end{lemma}
\begin{proof}
Set $U=\VA{A}$.
By the assumption, it suffices to show that $U_2\cdot T=0$.
Let $a\in A_\RR$ and $t\in T$.
First note that $(a|t)=(a\cdot u|t)=(u|a\cdot t)=0$ where $u=\omega/2$ is the unity of $B$. 
Hence we have $a_{(3)}t=(a|t)\1=0$.
We next consider the vector $a_{(0)}t\in V_3$.
We have $(a_{(1)}a)_{(1)}t+(a_{(2)}a)_{(0)}t
=a_{(2)}(a_{(0)}t)-a_{(1)}(a_{(1)}t)
+a_{(1)}(a_{(1)}t)-a_{(0)}(a_{(2)}t)$ 
 as a particular case of the Cauchy-Jacobi identity (\ref{CJ}).
Hence $a_{(2)}(a_{(0)}t)=0$, so $(a_{(0)}t|a_{(0)}t)=(t|a_{(2)}a_{(0)}t)=0$ and we have $a_{(0)}t=0$ by the positivity condition (I).
Hence $A_{(n)}T=0$ for all $n\geq 0$, which yields $U_{(n)}T$ for all $n\geq 0$ by the Cauchy-Jacobi identity again.
In particular, we have $U_2\cdot T=0$.
\end{proof}
\section{3-transposition Automorphisms}
We will present some results on the realization of 3-transposition automorphisms from a set of idempotents in the Griess algebra.
In particular, the notion of being $(\ccc,\cw)$-realizable by a VOA will be defined, though at this moment not many examples of such a case are known to the author except for the case with $\ccc=1/2$ and $\cw=1/2$.
\subsection{Fusion property of the Griess algebra}
\label{Fusion property of the Griess algebra}
Let $V$ be a VOA satisfying the condition (II) and let $B$ be the Griess algebra $V_2$.
We set
\begin{equation}
B^e=\{v\in B\,|\,(e|v)=0\}.
\end{equation}
If $c_e\ne 0$ then $B=B^e\oplus \CC e$.
For each complex number $h$, consider the space
\begin{equation}
B^e(h)=\{v\in B^e\,|\,e\cdot v=hv\}.
\end{equation}
This is the space of highest weight vectors with the conformal weight $h$ in $B$ since $L^e_2v=(e|v)\1=0$ on $B^e$.

In this situation, we have the following (cf.{} \cite{Mi1}).

\begin{lemma}
\sl
Let $V$ satisfy the conditions {\rm (I)} and {\rm (II)} and let $e\in B_\RR$ be a real idempotent.
Then we have
\begin{enumerate}[\rm (1)]
\item
The central charge $c_e$ is positive: $c_e>0$.
\item
$B^e$ is an orthogonal direct sum of eigenspaces $B^e(h)$ with $0\leq h\leq 2$.
\item
$L^e_{-1}B^e(0)=0$.
\end{enumerate}
\end{lemma}
\begin{proof}
The assertion (1) follows from (I) since $c_e=2(e|e)$ and $e\ne 0$.
Recall that $B=V_2$ is finite-dimensional by definition.
Since $(L^e_0 a|b)=(a|L^e_0 b)$ the action $L^e_0$ on $B_\RR$ is given by a symmetric matrix by taking an orthonormal basis, so it is diagonalizable with real eigenvalues for which the eigenspaces are orthogonal to each other.
Now suppose $B^e(h)\ne 0$.
Then there exists a nonzero real eigenvector $x\in B_\RR$ such that $L^e_0x=h x$.
Note that $(x|x)>0$ by the positivity (I) and by $L^e_1x=0$; The latter follows from the assumption (II).
Since  $2h(x|x)=(2L^e_0x|x)=([L^e_1,L^e_{-1}]x|x)=(L^e_1L^e_{-1}x|x)=(L^e_{-1}x|L^e_{-1}x)\geq 0$,  we have $h\geq 0$.
Now $\omega-e$ is also an idempotent and we have $L^{\omega-e}_0x=(2-h)x$.
Therefore $2-h\geq 0$ by the same argument as above.
Hence $0\leq h\leq 2$ and we have shown (2).
Finally if $h=0$ then we have $(L^e_{-1}x|L^e_{-1}x)=0$ for a real eigenvector $x\in B^e(0)$ so $L^e_{-1}x=0$ by (I) again, (cf.{} Lemma \ref{lemma: U2}).
\end{proof}

\vskip2ex

Next we wish to use the `fusion rules' in order to get more restrictions on the structure of $B$.
For convenience, we set
\newcommand{\even}{{\rm even}}
\newcommand{\odd}{{\rm odd}}
\begin{equation}
B_e[\barb{0}]=B^e(0)\oplus \CC e,\ B_e[\barb{1}]=B^e(\cw).
\end{equation}
Let us say that the idempotent $e$ has the {\it fusion property of binary type with respect to the spectrum $h=\cw$ in $B$} if the property
\begin{equation}
B_e[{\varepsilon_1}]\cdot B_e[{\varepsilon_2}]\subseteq B_e[{\varepsilon_1+\varepsilon_2}]
\label{fusion}
\end{equation}
is satsified for all $\varepsilon_1,\varepsilon_2\in \ZZ/2\ZZ$.

\vskip2ex

Note that $B^e(0)\cdot B^e(0)\subseteq B^e(0)$ and $B^e(0)\cdot B^e(\cw)\subseteq B^e(\cw)$, but $(B^e(\cw)\cdot B^e(\cw))\cap \CC e \ne \{0\}$ in general.
Indeed, for $x,y\in B^e(\cw)$, we have $(e|x\cdot y)=(e\cdot x|y)=\cw(x|y)\ne 0$.
\subsection{Pair of idempotents}
Let us fix real numbers $\ccc$ and $\cw$ such that $0<\ccc$ and $0<\cw<2$.
Let $e$ be an idempotent with the central charge $c=\ccc$ and let $x$ be a vector in the Griess algebra such that $x\in B^e(0)\oplus B^e(\cw)\oplus \CC$.
Let the corresponding decomposition be
\begin{equation}
x=x_e(0)+x_e(\cw)+x_e(2).
\end{equation}
Then we have
\begin{equation}
e\cdot x=\cw x_e(\cw)+\frac{4(e|x)}{\ccc}e,
\quad e\cdot(e\cdot x)=\cw^2x_e(\cw)+\frac{8(e|x)}{\ccc}e.
\end{equation}
Therefore, 
\begin{equation}
e\cdot(e\cdot x)=\cw e\cdot x+\frac{(8-4\cw )(e|x)}{\ccc}e.
\label{A.8}
\end{equation}

\vskip2ex

Let $e,f\in B$ be idempotents with the same central charge $c=\ccc$ and suppose that $e$ has the fusion property of binary type with respect to the spectrum $h=\cw$.
Suppose $f\in B^e(0)\oplus B^e(\cw)\oplus\CC e$ and $e\in B^f(0)\oplus B^f(\cw)\oplus \CC f$.
Since $f$ is an idempotent, we have
\begin{equation}
f\cdot f=2f.
\end{equation}
On the other hand, by (\ref{A.8}), we have
\begin{equation}
f\cdot(e\cdot f)=\cw e\cdot f+\frac{(8-4\cw)(e|f)}{\ccc}f.
\end{equation}
Consider the $B^e(\cw)$-components of these equalities and eliminate $f_e(0)\cdot f_e(\cw)$ from them.
Then we obtain
\begin{equation}
(1-\cw)\left(\frac{8(e|f)}{\ccc}-\cw\right)f_e(\cw)=0.
\end{equation}

If further $f$ has the fusion property of binary type with respect to the spectrum $h=\cw$, then we have the same equality with $e$ and $f$ being exchanged.
Thus we have the following result.
\begin{lemma}
\sl
Let the assumptions above are satisfied and suppose further that $\cw\ne 1$.
Then either $(e|f)=\cw \ccc/8$ or $e_f(\cw)=f_e(\cw)=0$.
\end{lemma}

Suppose $f_e(\cw)=0,\ e_f(\cw)=0$.
Then we have
\begin{equation}
e\cdot f=\frac{4(e|f)}{\ccc}e,\quad
e\cdot f=\frac{4(e|f)}{\ccc}f,
\end{equation}
Subtracting one from the other, we have
\begin{equation}
\frac{4(e|f)}{\ccc}(e-f)=0.
\end{equation}
If $e\ne f$ then $(e|f)=0$ hence $e\cdot f=0$.
Otherwise $e=f$.

Now consider the case when $(e|f)={\cw\ccc}/{8}$.
To describe the situation, we set
\begin{equation}
\sigma_e(x)=x-2x_e(\cw).
\end{equation}

\begin{lemma}
\sl
If $x\in B^e(0)\oplus B^e(\cw)\oplus \CC e$ then 
\begin{equation}
e\cdot x=\displaystyle\frac{\cw}{2}\left(x+\frac{8(e|x)}{\cw\ccc}e-\sigma_e(x)\right).
\end{equation}

\end{lemma}

\vskip2ex

\begin{lemma}
\sl
Let $e$ and $f$ be idempotents with the same central charge $c=\ccc$ such that $e\in B^f(0)\oplus B^f(\cw)\oplus\CC f$ and $f\in B^e(0)\oplus B^e(\cw)\oplus\CC e$.
Then 
\begin{equation}
\left(\frac{8}{\cw\ccc}(e|f)-1\right)(e-f)=\sigma_e(f)-\sigma_f(e).
\end{equation}
\end{lemma}

Therefore, if $e\ne f$ and $(e|f)=\cw c/8$ then $\sigma_e(f)=\sigma_f(e)$.
Note that in this case we have 
\begin{equation}
e\cdot f=\frac{\cw}{2}\left(e+f-\sigma_e(f)\right).
\end{equation}

\vskip2ex

The results obtained so far are summarized in the following theorem.

\begin{theorem}
\sl
Let $e,f\in B$ be indempotents with the same central charge $c=\ccc$ with the fusion property of binary type with respect to the spectrum $h=\cw\ne 1$.
Suppose $f\in B^e(0)\oplus B^e(\cw)\oplus \CC e$ and $e\in B^f(0)\oplus B^f(\cw)\oplus \CC f$.
Then one of the following holds:
\begin{enumerate}[\rm (1)]
\item
$(e|f)=\ccc/2$ and $e=f$.
\item
$(e|f)=0$ and $e\cdot f=0$.
\item
$(e|f)=\cw \ccc/8$ and $\sigma_e(f)=\sigma_f(e)$.
\end{enumerate}
\label{3cases}
\end{theorem}

\begin{note}
\rm
The results in this section in the case when $\ccc=1/2$ and $\cw=1/2$ and some of their consequences are due to Miyamoto \cite{Mi1}.
\end{note}
\subsection{$(\ccc,\cw)$-realizable groups}
\label{sect: realizable}
Let $e\in B_\RR$ be a real idempotent of the central charge $c=\ccc$ with the fusion property of binary type with respect to the spectrum $h=\cw$.
Suppose that $B=B^e(0)\oplus B^e(\cw)\oplus \CC e$. 
Then we may understand the operation $\sigma_e$ as a linear map
\begin{equation}
\sigma_e:B\rightarrow B,\ x\mapsto \sigma_e(x)=x-2x_e(\cw).
\end{equation}
The operator $\sigma_e$ acts by $(-1)^{\varepsilon}$ on $B_e[\varepsilon]$ for $\varepsilon\in\ZZ/2\ZZ$.

\begin{prop}
\sl
Let $e$ be an idempotent as above.
Then the map $\sigma_e$ is an automorphism of the algebra $B$ which preserves the bilinear form and we have $\tau\sigma_e\tau^{-1}=\sigma_{\tau(e)}$ for any automorphism $\tau$ of $B$.
\label{prop:autom conj}
\end{prop}
\begin{proof}
Suppose $x\in B_e[\varepsilon_1]$ and $y\in B_e[\varepsilon_2]$.
Then $x\cdot y\in B_e[{\varepsilon_1+\varepsilon_2}]$ by the fusion property.
Hence 
$\sigma_e(x\cdot y)=(-1)^{\varepsilon_1+\varepsilon_2}x\cdot y
=((-1)^{\varepsilon_1}x)\cdot ((-1)^{\varepsilon_2}y)
=(\sigma_e(x))\cdot (\sigma_e(y))$.
On the other hand, we have $(\sigma_e(x)|\sigma_e(y))=(-1)^{\varepsilon_1+\varepsilon_2}(x|y)=(x|y)$
since if $\varepsilon_1+\varepsilon_2=\barb{1}$ then $(x|y)=0$.
Now let $\tau$ be an automorphism of $B$ and let $x\in B_e(h)$ where $h=0,\cw$ or $2$.
Then since $\tau(e)\cdot \tau(x)=\tau(e\cdot x)=\tau(hx)=h\tau(x)$, we have $\tau(x)\in B_{\tau(e)}(h)$.
Therefore, $\tau\sigma_e\tau^{-1}$ acts by $(-1)^\varepsilon$ on the space $B_{\tau(e)}[\varepsilon]$.
This shows that $\tau\sigma_e\tau^{-1}=\sigma_{\tau(e)}$.
\end{proof}

\vskip2ex

Let $E_B$ be the set of the idempotents having all the properties assumed so far with respect to $c=\ccc$ and $h=\cw$.
For a pair $e,f\in E_B$, we write $e\sim f$ if $(e|f)=\cw\ccc/8$ and $e\perp f$ if $(e|f)=0$.
When $e\sim f$, we denote by $e\compo f$ the vector $\sigma_e(f)$, which is equal to $\sigma_f(e)$ by Theorem \ref{3cases}.

Let $D_B$ be the set of involutions of the form $\sigma_e$ for some $e\in E_B$ and let $G_B$ be the subgroup of $\Aut B$ generated by $D_B$ in $\Aut B$.

\begin{corollary}
\sl
The pair $(G_B,D_B)$ is a 3-transposition group.
\label{cor}
\end{corollary}
\begin{proof}
Let $e,f\in E_B$ be distinct idempotents.
If $e\perp f$ then since $f\in B_e[\barb{0}]$ we have $\sigma_e(f)=f$, hence $(\sigma_e\sigma_f)^2
=\sigma_e\sigma_f\sigma_e\sigma_f
=\sigma_{\sigma_e(f)}\sigma_f
=\sigma_f\sigma_f=1$.
If $e\sim f$ then we have $\sigma_e(f)=\sigma_f(e)$ by Theorem \ref{3cases}, hence $(\sigma_e\sigma_f)^3
=(\sigma_e\sigma_f\sigma_e)(\sigma_f\sigma_e\sigma_f)
=\sigma_{\sigma_e(f)}\sigma_{\sigma_f(e)}=1$.
\end{proof}

\vskip2ex

Now we wish to consider the action of $\sigma_e$ not only on the Griess algebra $B$ but also on the VOA $V$.
To this end, we need the fusion property regarding the whole space $V$.
Namely, we suppose that the space $V$ decomposes as $V=V_e[\barb{0}]\oplus V_e[\barb{1}]$ so that the following properties are satisfied:

\begin{enumerate}[(1)]
\item
$V_e[\varepsilon_0]_{(n)}V_e[\varepsilon_1]\subseteq V_e[\varepsilon_1+\varepsilon_2]$ for all $\varepsilon_1,\varepsilon_2\in\ZZ/2\ZZ$, $n\in\ZZ$.
\item
$V_e[\varepsilon]\cap V_2=B_e[\varepsilon]$ for $\varepsilon\in\ZZ/2\ZZ$.
\end{enumerate}

\noindent
Then the automorphism $\sigma_e$ of $B$ extends to an automorphism of the whole $V$ by the action which inverts the vectors in the component $V_e[\barb{1}]$.

Consider the set $E_V$ of idempotents satifsying these properties and assume further that

\begin{enumerate}[(1)]\setcounter{enumi}{2}
\item
$g(V_e[\varepsilon])=V_{g(e)}[\varepsilon]$ for all $g\in \Aut V$, $e\in E_V$ and $\varepsilon\in\ZZ/2\ZZ$.
\end{enumerate}

\noindent
Let $D_V$ be the set of involutions of the form $\sigma_e$ for some $e\in E_V$ and let $G_V$ be the subgroup of $\Aut V$ generated by $D_V$ in $\Aut V$.
Then by the same argument as in the proof of Proposition \ref{prop:autom conj} and Corollary \ref{cor} we can show the following.

\begin{prop}
\sl
Under the assumptions above, the map $\sigma_e$ is an automorphism of the VOA $V$ such that $\tau\sigma_e\tau^{-1}=\sigma_{\tau(e)}$ for any automorphism $\tau$ of $V$. Consequently the pair $(G_V,D_V)$ is a 3-transposition group.
\end{prop}

We will say that a 3-transposition group $(G,D)$ is $(\ccc,\cw)$-{\it realizable}\/ by a VOA $V$ if it is isomorphic to the 3-transposition subgroup $(G_V,D_V)$ of $\Aut V$ obtained as above.
We will often identify $(G_V,D_V)$ with $(G,D)$ without mentioning it explicitly.
We will denote $E_V$ simply by $E$ in the rest of the paper.
\subsection{Centerfreeness}
\label{Centerfreeness and indecomposability}
Let $(G,D)$ be a 3-transposition group which is $(\ccc,\cw)$-realizable by a  VOA $V$ satisfying (I) and (II).
Let $E$ be the set of associated idempotents defined as before:
\begin{equation}
E=\Bigl\{e\in B_\RR\,\Big|\,
\hbox{$e\cdot e=2e$, $2(e|e)=\ccc$, $V=V_e[\barb{0}]\oplus V_e[\barb{1}]$}\Bigr\}
\label{defE}
\end{equation}
so that $D=\{\sigma_e\in \Aut V\,|\,e\in E \}$.
The properties (1)--(3) in the preceding section are assumed.

\vskip2ex

Now since $E$ is stable under the action of an automorphism of $V$, we have a homomorphism $\Aut V\rightarrow \Symm E$ of groups where $\Symm E$ denotes the symmetric group on the set $E$.
Let $G_E$ be its image in $\Symm E$ and let $K$ be the kernel of $p|_{G}$: 
\begin{equation}
K=\Bigl\{\tau\in G\,\Big|\,\tau|_E=\id_E\Bigr\}.
\end{equation}
Thus we have an exact sequence of groups
\begin{equation}
1\rightarrow K\rightarrow G\rightarrow G_E\rightarrow 1.
\end{equation}
The following lemma is easily shown by a standard argument.
\begin{lemma}
\sl
The kernel $K$ is the center of $G$.
\end{lemma}
\begin{proof}
Let $\tau$ be an element of $K$.
Then, by the definition of $K$, it fixes the elements of $E$, so we have $\sigma_{\tau(e)}=\sigma_e$.
Since $\tau\sigma_e\tau^{-1}=\sigma_{\tau(e)}$, we see that $\tau$ commutes with all the elements in $D$. Hence $\tau$ is central in $G$ because $G$ is generated by $D$.
Conversely, let $\tau$ be central in $G$.
Let $e$ be an element of $E$ and set $f=\tau(e)$.
If there exists no element $x\in E$ such that $e\sim x$ then $D$ fixes $e$.
Since $D$ generates $G$ and $\tau\in G$, we have $e=f$ in this case.
So let $x\in E$ be such that $e\sim x$.
Then $\sigma_{f}(x)=\sigma_{\tau(e)}(x)=\tau\sigma_e\tau^{-1}(x)=\sigma_e(x)\ne x$, which implies $f\sim x$.
Therefore, $\sigma_x(e)=\sigma_e(x)=\sigma_{f}(x)=\sigma_x(f)$.
Applying $\sigma_x$, we have $e=f$.
Since $\tau(e)=e$ for an arbitrary $e$, we have $\tau\in K$.
\end{proof}

\vskip2ex

Since we are mainly interested in the classification of $(G,D)$ rather than of $V$, it will not be unnatural to assume the following condition.
\begin{description}
\item{(III)}
The VOA $V$ is generated by $B$ and the space $B$ spanned by $E$.
\end{description}

\begin{prop}
\sl
Let $V$ be a VOA satisfying the condition {\rm (I)--(III)}.
Then $G=G_E$ and the group $G$ is centerfree.
\end{prop}
\begin{proof}
If $V$ is generated by the set $E$ then the map $G\rightarrow G_E$ is an isomorphism since the action of an automorphism on $V$ is uniquely determined by its action on $E$.
Hence by the proposition above, the group $G$ is centerfree.
\end{proof}
\subsection{Indecomposability}
Now let $(G',D')$ and $(G'',D'')$ be 3-transposition groups.
Then their direct product is, by definition, the pair $(G,D)$ where
\begin{equation}
G=G'\times G'',\ D=(D'\times \{1\})\sqcup(\{1\}\times D''),
\end{equation}
which is again a 3-transposition group.
Let us call a 3-transposition group {\it indecomposable}\/ when it does not admit a nontrivial direct product decomposition.

\vskip2ex

The following lemma follows immediately from the construction of tensor product of VOAs.

\begin{prop}
\sl
If two 3-transposition groups are $(\ccc,\cw)$-realizable by VOAs $V'$ and $V''$ satisfying {\rm (I)--(III)}, respectively, then their direct product is $(\ccc,\cw)$-realizable by the tensor product $V'\otimes V''$ of VOAs.
\end{prop}

Conversely, suppose we have a nontrivial orthogonal decomposition of $E$, namely a decomposition $E=E'\sqcup E''$ into the disjoint union of nonempty subsets $E'$ and $E''$ such that $E'\perp E''$.
Let $V'$ and $V''$ be the the subVOA generated by $E'$ and $E''$ respectively.
By (I), $V'$ and $V''$ are simple VOAs, so is their tensor product $V'\otimes V''$ of VOAs.
Since $E'\perp E''$, we have $u'_{(n)}u''=0$ for all $u'\in V'$, $u''\in V''$ and $n\geq 0$.
Hence, by the property of the tensor product of VOAs, there is a homomorphism  $V'\otimes V''\rightarrow V$ of VOAs (cf.{} \cite{MN}), which is surjective by (III) and is an isomorphism by the simplicity.

Therefore, we may suppose without loss of generality that our group $(G,D)$ is indecomposable.
This means that the set $E$ does not have a nontrivial orthogonal decomposition.
We will say in such a case that the set $E$ is indecomposable by abuse of terminology.

On the other hand, if $E$ consists of one element then $V$ is merely isomorphic to $L(\ccc,0)$.
Hence we may further suppose $|E|\geq 2$ without loss of generality.
Then for any $e\in E$ there is at least one element $f\in E$ such that $e\sim f$, hence we have $\{e,f,\sigma_e(f)\}\subseteq E$; So in fact $|E|\geq 3$.

\vskip2ex

Thus we have come to the following assumption:
\begin{description}
\item{(IV)}
The set $E$ is indecomposable and $|E|\geq 3$.
\end{description}
Under this assumption, we have the following lemma.

\begin{lemma}
\sl
If $V$ satisfies {\rm (I)--(IV)} then the map $\sigma:E\rightarrow D$ is bijective.
\label{bijection}
\end{lemma}
\begin{proof}
It suffices to show the injectivity.
Suppose $\sigma_e=\sigma_f$.
By (IV) there is an element $g\in E$ such that $e\sim g$.
Then $g\ne \sigma_e(g)=\sigma_f(g)$ so $f\sim g$.
Hence $\sigma_g(e)=\sigma_e(g)=\sigma_f(g)=\sigma_g(f)$ showing $e=f$.
\end{proof}
\section{Algebras associated with 3-transposition Groups}
We will associate an algebra equipped a bilinear form to the Fischer space associated with any 3-transposition group so that it agrees, in case the group is $(\ccc,\cw)$-realizable by a VOA, with the Griess algebra of the VOA.
We will relate the positivity of the bilinear form to the condition on the least eigenvalue of the diagram of the Fischer space.
The realizability with respect to particular values of $\ccc$ and $\cw$ will be examined for a few examples of small groups.
\subsection{Fischer spaces}
\label{Fischer space}
We review the notion of the Fischer space associated with a 3-transposition group and its characterization due to Buekenhout.
See \cite{We1}, \cite{Ha1}, \cite{CH} and \cite{As} for detail.

A {\it partial linear space}\/ consists of a set $X$ called the set of points and a set $\Lc$ of subsets of $X$ called the set of lines such that any two points lie on at most one line and any line has at least two points.
Consequently for any lines $\ell_1$ and $\ell_2$ we have either $\ell_1\cap \ell_2=\emptyset$, $|\ell_1\cap \ell_2|=1$ or $\ell_1=\ell_2$.

Let $X$ be  partial linear space in which each line consists of three points.
Two distinct points $x$ and $y$ are called {\it collinear}\/ if they lie on a line.
The {\it diagram}\/ of a partial linear space $(X ,\Lc)$ is the collinearity graph of the space, namely the graph with the set of vertices $X$ such that two points are adjacent (connected by an edge) if and only if they are collinear. 
For a pair $x,y$ of distinct points of $X$, we write $x\sim y$ if $x$ and $y$ are collinear and $x\perp y$ if not.
When $x\sim y$, we denote by $x\compo y$ the third point on the line through $x$ and $y$.
\paragraph{Dual affine plane of order $2$}
\rm
Consider the set $X =\{x_{12},x_{13},x_{14},x_{23},x_{24},x_{34}\}$ of $6$ points and let the lines be given by 
\begin{equation}
\Lc=\{\ell_{1},\ell_2,\ell_3,\ell_4\},\quad \hbox{where}\quad \ell_i=\{x_{mn}\,|\,i\notin \{m,n\}\}.
\end{equation}
The space $(X ,\Lc)$ is called the dual affine plane of order $2$ since it is the dual (with points and lines being exchanged) to the affine plane $\FF_2^2$ of order $2$.
The associated diagram is formed by the vertices and the edges of an octahedron.

\paragraph{Affine plane of order $3$}
\rm
Consider the set $X =\{x_{ij}\,|\,0\leq i,j\leq 2\}$ of 9 points and let a $3$-set $\{x_{ij},x_{kl},x_{mn}\}$ be a line if and only if
$(i+k+m,j+l+n)\equiv (0,0)\ {\rm mod}\ 3$.
The resulting space with 9 points and 12 lines is the {\it affine plane of order $3$}.
The associated diagram is the complete graph on $9$ vertices.

\vskip2ex

According to Buekenhout (cf.{} \cite{CH}), a partial linear space $(X ,\Lc)$ for which the lines consist of three points is called an (abstract) {\it Fischer space}\/ if it satisfies the following propery:
\begin{description}
\item{(FS)}
For any two intersecting lines $\ell_1$ and $\ell_2$ the span of them is isomorphic either to the dual affine plane of order $2$ or to the affine plane of order $3$.
\end{description}

Now, let $(G,D)$ be a 3-transposition group.
The partial linear space associated with $(G,D)$ is the pair $(X ,\Lc)$ with $X =D$ the set of points and $\Lc$ the set of lines such that a subset $\ell\subseteq D$ is a line if and only if $\ell$ consists of three points which generate a subgroup of $G$ isomorphic to $\Sym_3$.
In other words, $\ell$ is a line if and only if $\ell$ is given by $\{\sigma,\tau,\tau^\sigma\}$ for some $\sigma,\tau\in D$ with the order of $\sigma\tau$ being $3$.
The following observation is fundamental (\cite{Fi}, cf.{} \cite{CH}, \cite{As}).

\begin{prop}[Fischer]
\sl
Let $(G,D)$ be a 3-transposition group.
Then the partial linear space associated with $(G,D)$ is a Fischer space.
\end{prop}

\vskip2ex

Now, conversely, let $(X ,\Lc)$ be an abstract Fischer space. Then, for $x\in X $, let $\sigma_x$ be the permutation of $X$ defined by setting
\begin{equation}
\sigma_x(y)
=\cases{
y&when $x$ and $y$ are not collinear,\cr
z&when $\{x,y,z\}$ is a line.\cr
}
\end{equation}
Then due to Bueknehout the pair $(G,D)$ where $D=\{\sigma_x\,|\,x\in X \}$ and $G=\langle D\rangle$ is a centerfree 3-transposition group.

The group constructed from the dual affine plane of order $2$ is isomorphic to the symmetric group $\Sym_4$ of degree $4$ and the one for the affine plane is of shape $3^2\split 2$.

\vskip4ex

\begin{remark}
\rm
If $(G,D)$ is $(\ccc,\cw)$-realizable then, thanks to Lemma \ref{bijection}, we may identify the Fischer space of $(G,D)$ as the space for which the set of points is the set $E$ (cf.{} (\ref{defE})) such that a 3-set $\ell$ is a line if and only if $\ell=\{e,f,e\compo f\}$ for some $e$ and $f$ with $e\sim f$ (recall that $e\compo f=\sigma_e(f)=\sigma_f(e)$).
Similarly, we can identify the diagram as the graph on the set $E$ with $e$ and $f$ being adjacent if and only if $e\sim f$.
\end{remark}
\subsection{Algebras associated with a Fischer space}
\label{Fischer space and Griess algebra}
A Fischer space is said to be {\it connected}\/ if the associated diagram is connected.

Let $X$ be a connected Fischer space.
Then the associated 3-transposition group acts transitively on the points.
We denote the number of the points of $X$ and the valency of the diagram as
\begin{equation}
\nu=|X|,\ k=|\{y\in X \,|\,y\sim x\}|.
\end{equation}
respectively.

Let $X$ be any Fischer space and let $\ccc$ and $\cw$ be any complex numbers.
Assign the symbol $\tildeb{x}$ to each $x\in X$ and regard the set $\{\tildeb{x}\,|\,x\in X \}$ as a basis of the vector space
\begin{equation}
\tildeb{B}=\bigoplus_{x\in X }\CC\tildeb{x}.
\end{equation}
We make this space into an algebra equipped with a bilinear form by Table \ref{table1}.
\begin{table}[t]
\begin{center}
\renewcommand{\arraystretch}{2}
\begin{tabular}{|c||c|c|}
\hline
cases&$\tildeb{x}\cdot \tildeb{y}$&$(\tildeb{x}|\tildeb{y})$\\
\hhline{|=#=|=|}
$x=y$&$2\tildeb{y}$&$\frac{\ccc}{2}$\\
\hline
$x\perp y$&$0$&$0$\\
\hline
$x\sim y$&$\frac{\cw}{2}(\tildeb{x}+\tildeb{y}-\widetildeb{x\compo y})$&$\frac{\cw\ccc}{8}$\\
\hline
\end{tabular}
\end{center}
\caption{Structure of the algebra $\tildeb{B}$}
\label{table1}
\end{table}

Obviously the multiplication is commutative and the bilinear form is symmetric.
It is less obvious but can be easily shown by using the characterizaton of the Fischer space (see Sect.{} \ref{Fischer space}) that the form is invariant with respect to the multiplication:
\begin{equation}
(\tildeb{e}\cdot \tildeb{f}|\tildeb{g})=(\tildeb{f}|\tildeb{e}\cdot \tildeb{g}).
\end{equation}

Finally, suppose $4+k\cw\ne 0$. Then, for any $y\in X $, we have
\begin{equation}
\biggl(\sum_{x\in X }\tildeb{x}\biggr)\cdot \tildeb{y}
=2\tildeb{y}+\frac{\cw}{2}\sum_{x\sim y}(\tildeb{x}+\tildeb{y}-\widetildeb{x\compo y})
=2\left(1+\frac{k\cw}{4}\right)\tildeb{y}.
\end{equation}
Therefore, the vector
\begin{equation}
\tildeb{\omega}=\frac{4}{4+k\cw}\sum_{x\in X }\tildeb{x}
\label{TildeOmega}
\end{equation}
satisfies $\tildeb{\omega}\cdot \tildeb{y}=2\tildeb{y}$ for any $y\in X $: The vector $\tildeb{\omega}$ is twice the unity of the algebra $\tildeb{B}$.

Now the value $c_{\tildeb{\omega}}=2(\tildeb{\omega}|\tildeb{\omega})$ is given by
\begin{equation}
\EQ{
2\left(\frac{4}{4+k\cw}\right)^2\sum_{x,y\in X}(\tildeb{x}|\tildeb{y})
&
=\frac{32}{(4+k\cw)^2}\left(\sum_{x\in X}(\tildeb{x}|\tildeb{x})+\sum_{x\in X}\smash{\sum_{\scriptstyle y\in X\atop\scriptstyle x\sim y}}(\tildeb{x}|\tildeb{y})\right)\vphantom{\sum_{\scriptstyle y\in X\atop\scriptstyle x\sim y}}\cr
&
=\frac{32}{(4+k\cw)^2}\left(\frac{\ccc}{2}\nu+\frac{\cw\ccc}{8}\nu k\right)
=\frac{4\ccc \nu}{4+k\cw}\cr
}
\label{rank}
\end{equation}
which, if the associated 3-transposition group is $(\ccc,\cw)$-realizable by a VOA, turns out to be the rank of the VOA.

To summarize, we have obtained the following result.

\begin{prop}
\sl
For any Fischer space $X $, the algebra $\tildeb{B}$ is a commutative nonassociative algebra equipped with a symmetric invariant bilinear form.
If $4+k\cw\ne 0$ then the algebra has a unity.
\end{prop}

Now let $(G,D)$ be a 3-transposition group which is $(\ccc,\cw)$-realizable by a VOA satisfying (I)--(IV).
We identify the set $E$ of idempotents given by (\ref{defE}) with the Fischer space $X$ associated with $(G,D)$.
Consider the map
\begin{equation}
\pi:\tildeb{B}\rightarrow B,\quad \tildeb{e}\mapsto e.
\end{equation}
Then obviously
\begin{prop}
\sl
Suppose that a 3-transposition group $(G,D)$ is $(\ccc,\cw)$-realizable by a VOA $V$.
Then the map $\pi$ defined above is a surjective homomorphism of algebras that perserves the bilinear form and the kernel of $\pi$ agrees with the radical of the form $(\ |\ )$.
Furthermore, the rank of the VOA $V$ is equal to $c^{}_V={4\nu \ccc}/({4+k\cw})$.
\label{covers}
\end{prop}

\begin{remark}
\rm
The structure of an algebra with a bilinear form such as in the theorem is unique up to the choice of parameters $\ccc$ and $\cw$.
\end{remark}

\begin{note}
\rm
Our algebra $\tildeb{B}$ or $B$ differs from those considered by Norton \cite{No} in the existence of a unity.
Our algebra is a generalization of those considered in Dong et al.{} \cite{DLMN} (see Sect.{} \ref{root systems}).
\end{note}

Thus the multiplication table of the Griess algebra is determined by the Fischer space. 
On the other hand, the bilinear form on $B$ is determined by the diagram of the 3-transposition group.
Indeed, the Gram matrix on the space $\tildeb{B}$ is the matrix 
\begin{equation}
\frac{\ccc}{2}\left(I+\frac{\cw}{4}A\right)
\end{equation}
where $A$ is the adjacency matrix of the diagram.
Let $s$ be the least eigenvalue of the matrix $A$.

Since the form $(\ |\ )|_{B_\RR}$ is positive-definite, we have the following: 

\begin{theorem}
\sl
Suppose that a 3-transposition group $(G,D)$ is $(\ccc,\cw)$-realizable by a VOA satisfying (I)--(IV).
Then the least eigenvalue $s$ of the matrix $A$ satisfies $s\geq -4/\cw$, and the kernel of the map $\pi$ is nontrivial if and only if $s=-4/\cw$.
\label{least eigenvalue}
\end{theorem}
\subsection{Examples for small groups}
\label{Examples}
We will use the symbol $p^m$ to denote the elementary abelian group $(\ZZ/p\ZZ)^n$ as frequently used in finite group thoery.
For a group $G$ and a $G$-module $N$, the symbol $N\split G$ means a split extension of $G$ by $N$.
\subsubsection*{Projective line of order $2$ (The group $\Sym_3$)}
\label{Projective line of order $2$}
Let $G=\Sym_3$ and let $D$ be the set of transpositions:
\begin{equation}
G=\langle D\rangle,\ D=\{(12),(13),(23)\}.
\end{equation}
Then $(G,D)$ is an indecomposable centerfree 3-transposition group.
The associated Fischer space is $X =\{e_{1},e_{2},e_{3}\}$ with only one line: $\Lc=\{X \}$.
The associated algebra is $B=\CC e_{1}\oplus \CC e_{2} \oplus \CC e_{3}$ with the multiplication table
\begin{equation}
e_{i}\cdot e_{j}=\cases{\hfill 2e_j\hfill&if $i=j$,\cr\displaystyle \frac{\cw}{2}(e_i+e_j-e_k)&if $i\ne j$ and $\{i,j,k\}=\{1,2,3\}$.
\cr}
\end{equation}

Suppose that $(G,D)$ is $(\ccc,\cw)$-realizable by a VOA $V$.
Then
\begin{equation}
\omega=\frac{4}{4+2\cw}(e_1+e_2+e_3)
\end{equation}
is the conformal vector and the rank is given by
\begin{equation}
c^{}_V=\frac{12\ccc}{4+2\cw}=\frac{6\ccc}{2+\cw}.
\end{equation}
Set $\xi=e_1$ and $\eta=\omega-e_1$.
Then they are idempotents with the central charge
\begin{equation}
c_\xi=\ccc\hbox{\quad and\quad}{c_\eta}=\frac{6\ccc}{2+\cw}-\ccc=\frac{(4-\cw)\ccc}{2+\cw},
\end{equation}
respectively.
We have $\VA{\xi}\simeq L(c^{}_\xi,0)$ and $\VA{{\eta}}\simeq L({c_\eta},0)$.
Since the actions of corresponding Virasoro algebras commute with each other, 
we have an embedding of VOAs
\begin{equation}
L({c^{}_\xi},0)\otimes L({c_\eta},0)\hookrightarrow V
\end{equation}
onto a full subVOA.

\begin{prop}
\sl
The group $\Sym_3$ is $(1/2,1/2)$-realizable.
\end{prop}

For $\ccc=1/2$ and $\cw=1/2$, we have $c_\xi=1/2$ and $c_\eta=7/10$.
By inspecting the decomposition with respect to $L(1/2,0)\otimes L(7/10,0)$  it is not difficult to see that our algebra $\tildeb{B}$ is isomorphic to the Griess algebra of
\begin{equation}
V=\biggl(L(1/2,0)\otimes L(7/10,0)\biggr)\oplus\biggl(L(1/2,1/2)\otimes L(7/10,3/2)\biggr).
\end{equation}

The existence of a VOA structure on this space is easily seen by looking at the VOA $V^+_{\sqrt{2}A_2}$; See the next subsection and Sect.{} \ref{root systems}.
The VOA structure as well as its representation theory is studied by Lam and Yamada \cite{LY}.

\vskip2ex

Now let us turn to the case when $\ccc=1/2$ and $\cw=1/16$. 
A recent result of Miyamoto \cite{Mi2} guarantees the following.

\begin{prop}
\sl
The group $\Sym_3$ is $(1/2,1/16)$-realizable.
\end{prop}

The rank of the VOA is given by
\begin{equation}
c^{}_V=\frac{4\cdot 3\cdot \ccc}{4+2\cdot \cw}
=\frac{16}{11}=\frac{1}{2}+\frac{21}{22}
\end{equation}
and the VOA contains a subspace
\begin{equation}
\biggl(L(1/2,0)\otimes L(21/22,0)\biggr)\oplus\biggl(L(1/2,1/16)\otimes L(21/22,31/16)\biggr).
\end{equation}

\begin{note}
\rm
By \cite{Mi2}, such a VOA is constructed in the moonshine module $V^\natural$  whose automorphism group is the Monster, the largest sporadic finite simple group.
Indeed, by Conway \cite{Co}, an idempotent called the transposition axis is assigned to each involution in the conjugacy class 2A in the Monster which turns out to be an idempotent with the central charge $c=1/2$ (\cite{Mi1}).
Now take two elements of 2A such that their product falls into the class 3C then, by \cite{Mi2}, subalgebra generated by these idempotents must have the properties for the Griess algebra which realize $\Sym_3$ for $c=1/2$ and $h=1/16$ as above.
Recently C.-H.{} Lam \cite{La} has given a construction of the same VOA independent of the Monster.
The multiplication rule $e\cdot f=(1/4)(e+f-e\compo f)$, up to normalization, is already given in the row of the class 3C in Table 3 of \cite{Co}.
\label{Miyamoto 3C}
\end{note}

\begin{remark}
\rm
It is more or less obvious from the consideration above that the group $\Sym_3$ is both $(7/10,3/2)$-realizable and $(21/22,31/16)$-realizable.
\end{remark}
\subsubsection*{Dual affine plane of order $2$ (The group $\Sym_4$)}
\label{Dual affine plane}
Let $G=\Sym_4$ and let $D$ be the set of transpositions:
\begin{equation}
G=\langle D\rangle,\ D=\{(12),(13),(14),(23),(24),(34)\}.
\end{equation}
Then $(G,D)$ is an indecomposable centerfree 3-transposition group.
It can also be described as $2^2\split\Sym_3$.
The associated Fischer space is given in Sect.{} \ref{Fischer space}.
We have $\nu=6$ and $k=4$ for the diagram assocaited with the Fischer space.

Suppose that $(G,D)$ is $(\ccc,\cw)$-realizable by a VOA $V$.
Then
\begin{equation}
\omega=\frac{4}{4+4\cw}(e_{12}+e_{13}+e_{14}+e_{23}+e_{24}+e_{34})
\end{equation}
is the conformal vector and the rank is given by
\begin{equation}
c^{}_V=\frac{24\ccc}{4+4\cw}=\frac{6\ccc}{1+\cw}.
\end{equation}
Consider the subalgebra of the Griess algebra spanned by the $e_{12},e_{23},e_{13}$.
Then it is the algebra associated with the projective line of order $2$ as described in the preceding subsection.
Let $e_{123}$ be the conformal vector of this subalgebra and set $e_{1234}=\omega$.
Set ${\xi}=e_{12}$, ${\eta}=e_{123}-e_{12}$ and ${\zeta}=e_{1234}-e_{123}$.

Then the three vectors ${\xi},{\eta},{\zeta}$ generate mutually commutative Virasoro actions.
Hence we have an embedding $L(c_{\xi},0)\otimes L({c_\eta},0)\otimes L({c_\zeta},0)\hookrightarrow V$, where the central charge is given by
\begin{equation}
c_{\xi}=\ccc,\ {c_\eta}=\frac{(4-\cw)\ccc}{2+\cw},\ {c_\zeta}=\frac{6\ccc}{(1+\cw)(2+\cw)}
\end{equation}
respectively.

\begin{prop}
\sl
The group $\Sym_4$ is $(1/2,1/2)$-realizable.
\end{prop}

The group is realized by the VOA $V^+_{\sqrt{2}{A_2}}$; See Sect.{} \ref{root systems}.

\begin{note}
\rm
The author does not know whether $\Sym_4$ is $(1/2,1/16)$-realizable or not.
The consideration as in the proof above does not prohibit the existence of such VOA numerically.
\end{note}

\begin{remark}
\rm
There exists an eigenvector with nonzero norm with respect to $\eta$ whose eigenvalue is equal to ${\cw(1-\cw)}/{(2+\cw)}$, which  has to be positive by the assumption (I). 
Therefore, if $\Sym_4$ is $(\ccc,\cw)$-realizable then $0<\cw<1$.
\label{less than 1}
\end{remark}
\subsubsection*{Affine plane of order $3$ (The group of shape $3^2\split 2$)}
Let $X$ be the affine plane of order $3$ and let $(G,D)$ be the associated 3-transposition group.
Then the group $G$ has the structure $3^2\split2$, where the generator of $2=\ZZ/2\ZZ$ acts by inverting the vectors of $3^2=\FF_3^2$.
The structure of the Fischer space is described in Sect.{} \ref{Fischer space}.
We have $\nu=9$ and $k=8$ for the associated diagram.

Suppose that $(G,D)$ is $(\ccc,\cw)$-realizable by a VOA $V$.
Then
\begin{equation}
\omega=\frac{4}{4+8\cw}\sum_{0\leq i,j\leq 2}e_{ij}
\end{equation}
is the conformal vector and the rank is given by
\begin{equation}
c^{}_V=\frac{4\cdot 9\cdot \ccc}{4+8\cdot \cw}=\frac{9\ccc}{1+2\cw}.
\end{equation}
Consider the subalgebra of the Griess algebra $B_{0*}$ spanned by the $e_{00},e_{01},e_{02}$.
Then $B_{0*}$ is isomorphic to the algebra associated with the projective line of order $2$ as described in the preceding subsection.
Let $e_{0*}$ be the conformal vector of this subalgebra and $e_{**}$ be the conformal vector $\omega$.

Set $\xi=e_{00}$, ${\eta}=e_{0*}-e_{00}$ and $\zeta=e_{**}-e_{0*}$.
Then the three vectors $\xi,{\eta},{\zeta}$ give mutually commutative Virasoro actions.
Hence we have an embedding
\begin{equation}
L(c^{}_\xi,0)\otimes L({c_\eta},0)\otimes L({c_\zeta},0)\hookrightarrow V,
\end{equation}
where the central charge is given by
\begin{equation}
c^{}_\xi=\ccc,\ {c_\eta}=\frac{(4-\cw)\ccc}{2+\cw},\ {c_\zeta}=\frac{3(4-\cw)\ccc}{(1+2\cw)(2+\cw)}
\end{equation}
respectively.

\begin{prop}
\sl
The group $3^2\split 2$ is {\it not}\/ $(1/2,1/2)$-realizable.
\label{prop: NotRealizable}
\end{prop}
\begin{proof}
Consider the vectors
\begin{equation}
{\eta}=\frac{4}{5}(e_{00}+e_{01}+e_{02})-e_{00},\quad
w=e_{10}+e_{11}+e_{12}-e_{20}-e_{21}-e_{22}.
\end{equation}
Then, as we have already observed, the vector ${\eta}$ is an idempotent which generates a subVOA isomorphic to $L(7/10,0)$.
Further, we have 
${\eta}\cdot w=({7}/{10})w$ and $(w|w)=3/2$.
The former means $L^{{\eta}}_0w=(7/10)w$ and the latter implies $w\ne 0$.
Hence the vector $w$ is a highest weight vector of conformal weight $7/10$ with respect to the action of the Virasoro algebra with the central charge $c=7/10$.
However, this is impossible because the conformal weight of unitary irreducible highest weight representations of the central charge $c=7/10$ are only $0$, $1/10$, $3/5$, $3/2$, $3/80$ and $7/16$ (\cite{FQS}).
\end{proof}

\vskip2ex

\begin{note}
\rm
The author does not know whether $3^2\split 2$ is $(1/2,1/16)$-realizable or not.
In fact, the eigenvalues as in the proof are $0,31/16,1/11,5/176,21/176$, which are all allowed values of the conformal weight of unitary highest weight representations of $L(21/22,0)$.
\end{note}

\begin{remark}
\rm
By the same reason as Remark \ref{less than 1}, we see that if $3^2\split 2$ is $(\ccc,\cw)$-realizable then $0<\cw<1$.
\end{remark}
\section{Classification of $(1/2,1/2)$-realizable Groups}
We will restrict our attention to the case $\ccc=1/2$ and $\cw=1/2$ and show the main result.
\subsection{3-transposition groups of symplectic type}
\label{3-transposition group of symplectic type}
Let $(G,D)$ be a centerfree 3-transposition group.
It is called {\it of symplectic type}\/ if the affine plane of order $3$ does not occur in the associated Fischer space.
Such groups are classified by J.I.{} Hall in \cite{Ha1} and \cite{Ha2}.

\vskip2ex

\begin{theorem}[J.I.{} Hall]
\sl
An indecomposable centerfree $3$-transposition group of symplectic type is isomorphic to the extension of one of the groups {\rm $\Sym_3$; $\Sym_n$, ($n\geq 5$); $\Sp_{2n}(2)$, $(n\geq 3$); $\Orth^+_{2n}(2)$, ($n\geq 4$);} and {\rm $\Orth^-_{2n}(2)$, ($n\geq 3)$,} by the direct sum of copies of the natural module.
\label{theorem: hall}
\end{theorem}

\vskip2ex

Here the natural module, which we will denote by $\EE$ in the sequel, is isomorphic to $2^{2n}$ for $\Orth^{\pm}_{2n}(2)$, or $\Sp_{2n}(2)$.
Since $\Sym_{2n+1}$ and $\Sym_{2n+2}$ are embedded in the symplectic transformation on the space $2^{2n}$, we understand that the natural module for these groups is $2^{2n}$.
Note that $\Sym_4\simeq 2^2\split\Sym_3$.
We denote by $G^*$ the group without an extension by the natural modules.

\vskip2ex

Let us first consider the case without an extension.
Then the group $G=G^*$ is a rank 3 permutation group on the set $D$ and the associated diagram is a strongly regular graph (cf.{} \cite{Fi},\cite{We1}).
Let $D^{\Point}$ be the elements of $D$ which commutes with a fixed element of $D$ and $D^{\Line}$ be that with two fixed noncommuting elements of $D$.
(They are usually denoted by $D_d$ and $D_{d,e}$ respectively.)
The elements of $D^{\Line}$ actually commute with 3 elements that correspond to a line in the associated Fischer space.
We set $G^{\Point}=\langle D^{\Point}\rangle$ and $G^{\Line}=\langle D^{\Line}\rangle$.
Table \ref{table2} summarizes the inductive structure of those 3-transposition groups of symplectic type (cf.{} \cite{We2}).

{\small
\begin{table}[t]
\begin{center}
\renewcommand{\arraystretch}{1.5}
\begin{tabular}{|rc||c|c|c|}
\hline
$G$&&$|D|$&$G^{\Point}$&$G^{\Line}$\\
\hline
$\Sym_3$&&3&---&---\\
\hline
$\Sym_4$&&6&$S_2$&---\\
\hline
$\Sym_n$&\hskip-.5em$n\geq 5$&$\frac{n(n-1)}{2}$&$\Sym_{n-2}$&$\Sym_{n-3}$\\
\hline
$\Orth^+_{2n}(2)$&\hskip-.5em$n\geq 4$&$2^{2n-1}-2^{n-1}$&$\Sp_{2n-2}(2)$&$\Orth^-_{2n-2}(2)$\\
\hline
$\Orth^-_{2n}(2)$&\hskip-.5em$n\geq 3$&$2^{2n-1}+2^{n-1}$&$\Sp_{2n-2}(2)$&$\Orth^+_{2n-2}(2)$\\
\hline
$\Sp_{2n}(2)$&\hskip-.5em$n\geq 3$&$2^{2n}-1$&$2^{2n-2}\split \Sp_{2n-2}(2)$&$\Sp_{2n-2}(2)$\\
\hline
\end{tabular}
\end{center}
\caption{Inductive structure of 3-transposition groups of symplectic type}
\label{table2}
\end{table}
}
From Table \ref{table2}, we can compute the standard parameters $(\nu,k,\lambda,\mu)$ of the associated diagram by standard techniques as follows (cf.{} \cite{Bos},\cite{Hi}).
We set
\begin{equation}
\nu=|D|,\ \barb{k}=|D^{\Point}|,\ \barb{\mu}=|D^{\Line}|.
\end{equation}
Then we have
\begin{equation}
k=\nu-\barb{k}-1,\ 
\lambda=\barb{\mu}+2k-\nu,\ \mu=\frac{k(k-1-\lambda)}{\barb{k}}.
\end{equation}
The eigenvalues of the diagram are $k$ and the roots of
\begin{equation}
x^2-(\lambda-\mu)x+(\mu-k)=0.
\end{equation}
Let $r$ be the greater root and $s$ the lesser one:
\begin{equation}
r=\frac{(\lambda-\mu)+\sqrt{(\lambda-\mu)^2-4(\mu-k)}}{2},\ 
s=\frac{(\lambda-\mu)-\sqrt{(\lambda-\mu)^2-4(\mu-k)}}{2}.
\end{equation}
For our later purpose, we also consider the multiplicity $g$ of the least eigenvalue $s$; It is given by the following formula:
\begin{equation}
g=\frac{k+r(\nu-1)}{r-s}.
\end{equation}

Now, let us turn to the case with an extension: By the theorem of Hall mentioned above, we have $G=\EE^m\split G^*$ for some $m$ where $\EE$ is the natural module over $G^*$.

Let $X^*$ be the Fischer space associated with the group $G^*$ and let $X$ be the set $2^m\times X^*$.
Define the structure of a partial linear space on the set $X$ by letting a 3-set $\{(p,x),(q,y),(r,z)\}$ be a line if and only if $\{x,y,z\}$ is a line of $X^*$ and $p+q+r\equiv 0$ modulo $2$.
It is easy to check that the space $X=2^m\times X^*$ is a Fischer space of symplectic type.
By the theorem of Hall, we have the following.

\begin{prop}
\sl
The partial linear space $X=2^m\times X^*$ is a Fischer space of symplectic type isomorphic to the Fischer space associated with the group $G=\EE^m\split G^*$.
\label{covering}
\end{prop}

In particular, the diagram associated with $G$ or $X$ is identified with the graph on the set $2^m\times X^*$ with two vertices $(p,x)$ and $(q,y)$ being adjacent if and only if $x\sim y$.
Hence we have the following.

\begin{lemma}
\sl
The Gram matrix of the diagram of the Fischer space associated with $G=\EE^m\split G^*$ is given by
\begin{equation}
A=\underbrace{\pmatrix{
1&1\cr
1&1\cr
}\otimes\cdots\otimes\pmatrix{
1&1\cr
1&1\cr
}}_{\textstyle m}\otimes A^*,
\end{equation}
where $A^*$ is the adjacency matrix of the diagram associated with $X^*$.
\label{tensor}
\end{lemma}

The eigenvalues of $A$ are simply $0$ and $2^m$ times the eigenvalues of $A^*$.
In particular, we have $s=2^m s^*,\ g=g^*$.

The results are summarized in Table \ref{table3}.

{\small
\begin{table}[t]
\begin{center}
\renewcommand{\arraystretch}{1.5}
\begin{tabular}{|rc||c|c|c|c|c|}
\hline
$G$&&$\nu$&$k$&$\lambda$&$s$&$g$\\
\hhline{|==#=|=|=|=|=|}
$\Sym_3$&&$3$&$2$&$1$&$-1$&$2$\\
\hline
$\Sym_n$&\hskip-.5em$n\geq 4$&$\frac{n(n-1)}{2}$&$2n-4$&$n-2$&$-2$&$\frac{n(n-3)}{2}$\\
\hhline{|==#=|=|=|=|=|}
$\Orth^+_{2n}(2)$&\hskip-.5em$n\geq 4$&$\scriptstyle2^{2n-1}-2^{n-1}$&$\scriptstyle2^{2n-2}-2^{n-1}$&$\scriptstyle2^{2n-3}-2^{n-2}$&$\scriptstyle-2^{n-2}$&$\frac{2^{2n}-4}{3}$\\
\hline
$\Orth^-_{2n}(2)$&\hskip-.5em$n\geq 3$&$\scriptstyle2^{2n-1}+2^{n-1}$&$\scriptstyle2^{2n-2}+2^{n-1}$&$\scriptstyle2^{2n-3}+2^{n-1}$&$\scriptstyle-2^{n-1}$&$\frac{2^{2n}+3\cdot 2^{n}+2}{6}$\\
\hline
$\Sp_{2n}(2)$&\hskip-.5em$n\geq 3$&$\scriptstyle2^{2n}-1$&$\scriptstyle2^{2n-1}$&$\scriptstyle2^{2n-2}$&$\scriptstyle-2^{n-1}$&$\scriptstyle2^{2n-1}+2^{n-1}-1$\\
\hhline{|==#=|=|=|=|=|}
$\EE^m\split G^*\hskip-1.7em$&&$2^m\nu^*$&$2^mk^*$&$2^m\lambda^*$&$2^ms^*$&$g^*$\\
\hline
\end{tabular}
\end{center}
\caption{Parameters of the diagram}
\label{table3}
\end{table}
}
\subsection{Groups related to root systems}
\label{root systems}
Let $R$ be a root system of simply-laced type.
We denote the root lattice by the same symbol $R$.
Let $\sqrt{2}R$ be the root lattice with the norm being multiplied by $2$, which is denoted by $R(2)$ in other areas of mathematics.
Let $V_{\sqrt{2}R}$ denote the VOA associated with this lattice and let $V=V^+_{\sqrt{2}R}$ be the fixed-point subspace with respect to the involution $\theta$ which is a lift of the $-1$ isometry of the lattice (\cite{FLM}).

It is well known that $V^+_{\sqrt{2}R}$ is a full subVOA of $V_{\sqrt{2}R}$ satisfying the properties (I)--(IV) for $\ccc=1/2$ and $\cw=1/2$.
In particular, for an idempotent in $E$, we have
\begin{equation}
V=V(0)\oplus V(1/2)
\end{equation}
where $V(0)$ and $V(1/2)$ are the sum of components isomorphic to $L(1/2,0)$ and $L(1/2,1/2)$ respectively as a module over the Virasoro algebra with the central charge $1/2$ generated by the idempotent.
The properties (1)--(3) in Sect.{} \ref{sect: realizable} are satisfied. (The property (1) is a consequence of the fusion rules of such modules.)

The structure of the Griess algebra of $V^+_{\sqrt{2}R}$ is described by Dong et al.{} \cite{DLMN}.
Let us first recall this in our terms.

Let $R^+$ be the set of positive roots and let $X^*$ be the partial linear space on the set $R^+$ such that a 3-set is a line if and only if it is the set of the three positive roots in a subsystem of type $A_2$.
Now apply the construction in Proposition \ref{covering} to $X^*$ with $m=1$.
Then the algebra $\tildeb{B}$ associated with the resulting Fischer space $X=2\times X^*$ covers the Griess algebra of $V^+_{\sqrt{2}R}$ when $\ccc=1/2$ and $\cw=1/2$ (see Proposition \ref{covers}).
The algebra $\tildeb{B}$ is nothing else but the algebra $B^+$ in \cite{DLMN}.

\vskip2ex

Let us denote by $B$ the Griess algebra of $V^+_{\sqrt{2}R}$ and let $B^*$ be the subalgebra generated by the image of $X^*$.
The spaces $X$ and $X^*$ are the Fischer space associated with the groups $G$ and $G^*$ as given in Table \ref{Table: root}.
Consider the subVOA $V^*$ generated by the subalgebra $B^*$.
We will show below that $V^*_2$ agrees with $B^*$.

Let $\tildeb{\omega}$ be twice the unit of the larger algebra $\tildeb{B}$ and let $\tildeb{\omega}^*$ be that of $\tildeb{B}^*$.
Set
\begin{equation}
\tildeb{\eta}=\tildeb{\omega}-\tildeb{\omega}^*.
\end{equation}
Consider the eigenspace decomposition of the large algebra $\tildeb{B}$ with respect to the adjoint action by $\tildeb{\eta}$.
Obviously, the small algebra $\tildeb{B}^*$ is contained in the eigenspace with the eigenvalue $0$.
Then the induced action of $\tildeb{\eta}$ on the quotient space $\tildeb{B}/\tildeb{B}^*$ is given by the following matrix:
\begin{equation}
\left(2-\frac{k}{4+k\cw}\right)I+\frac{1}{4+k\cw}A^*.
\end{equation}
where $I$ is the identity matrix of size $|X^*|$ and $A^*$ is the adjacency matrix of the diagram associated with $X^*$.
Therefore, when $\cw=1/2$, this turns out to 
\begin{equation}
\frac{2}{8+k}\left(8I+A^*\right).
\label{matrix on quotient space}
\end{equation}
This matrix is positive since the least eigenvalue of the matrix $A^*$ is greater than or equal to $-4$.
Hence we have proved the following.
\begin{lemma}
\sl
The subspace ${B}^*$ agrees with the eigenspace with the eigenvalue $0$ with respect to the adjoint action by $\eta$.
\end{lemma}

{\small
\begin{table}
\begin{center}
\renewcommand{\arraystretch}{2}
\begin{tabular}{|c||c|c|c|c|}
\hline
$R$&$G^*$&$G$&$c_{\eta}$&eigenvalues of $\eta$ on $B$\\
\hhline{|=#=|=|=|=|}
$A_{n-1}$&$\Sym_n$&$\EE\split \Sym_n$&$\frac{2(n-1)}{n+2}$&$0,\ \frac{6}{n+2},\frac{n+4}{n+2},2$\\
\hline
$D_{n}$&$\EE\split\Sym_n$&$\EE^2\split \Sym_n$&$1$&$0,\ \frac{4}{n},\ 1,\ 2$\\
\hhline{|=#=|=|=|=|}
$E_{6}$&$\Orth^-_6(2)$&$2^6\split \Orth^-_6(2)$&$\frac{6}{7}$&$0,\ \frac{7}{5},\ 2$\\
\hline
$E_{7}$&$\Sp_6(2)$&$2^6\split \Sp^-_6(2)$&$\frac{7}{10}$&$0,\ \frac{3}{5},\ 2$\\
\hline
$E_{8}$&$\Orth^+_8(2)$&$2^8\split \Orth^+_8(2)$&$\frac{1}{2}$&$0,\ \frac{1}{2},\ 2$\\
\hline
\end{tabular}
\end{center}
\caption{Groups related to root systems}
\label{Table: root}
\end{table}
}
Then by Lemma \ref{lemma: U2}, we see the following.

\begin{lemma}
\sl
The Griess algebra of the subVOA $V^*$ generated by $B^*$ agrees with the algebra $B^*$.
\end{lemma}

\vskip2ex

Therefore, we have verified the following:

\begin{prop}
\sl
The groups $G$ and $G^*$ in Table \ref{Table: root} are  $(1/2,1/2)$-realizable by a VOA satisfying {\rm (I)}--{\rm (IV)}.
\label{prop: root}
\end{prop}

There are some more groups which are $(1/2,1/2)$-realizable.
The complete list will be given in the next section.

\begin{note}
\rm
For $E_8$ type root system, the central charge of $\eta$ is equal to $1/2$.
This implies that for $E_8$ we have more idempotent with the central charge $1/2$.
The detailed account for the Griess algebra of $V^+_{\sqrt{2}E_8}$ concluding the isomorphism ${\rm Aut}\,V^+_{\sqrt{2}E_8}\simeq \Orth^+_{10}(2)$ is given in \cite{Gr}. 
The full automorphism group of $V^+_{\sqrt{2}R}$ for other types of root system is described in \cite{Sh}.
\label{note: Griess}
\end{note}

\begin{note}
\rm
For the $DE$ type root system, we have the following formula:
\begin{equation}
g=|R_+|-\frac{\ell(\ell+1)}{2},
\end{equation}
where $R$ is the root system, $R_+$ a set of positive roots, $\ell$ the rank and $g$ the multiplicity of the least eigenvalue of the graph defined on $R$ by setting $\alpha\sim\beta$ if and only if $\langle\alpha,\beta\rangle=\pm 1$.
\end{note}
\subsection{Classification of $(1/2,1/2)$-realizable groups}
\label{Classification for $1/2$ and $1/2$}
Let $(G,D)$ be a centerfree 3-transposition group which is $(1/2,1/2)$-realizable.

Recall Proposition \ref{prop: NotRealizable}.
This implies that the associated Fischer space is of symplectic type.
On the other hand, Theorem \ref{least eigenvalue} says that the least eigenvalue of the adjacency matrix $A$ of the diagram of $X$ must be greater than or equal to $-8$.

Among the groups in Table \ref{table3}, the cases with $s\geq -8$ are only:
\begin{equation}
\Sym_3,\ \Sym_n,\ (n\geq 5),\ \Orth^+_{8}(2),\ \Orth^+_{10}(2),\ \Orth^-_{6}(2),\ \Orth^-_{8}(2),\ \Sp_{6}(2),\ \Sp_{8}(2).
\end{equation}
so the allowed extensions are as follows.
\begin{equation}
\Sym_4=2^2\split \Sym_3,\ \EE\split \Sym_n,\ \EE^2\split \Sym_n,\ (n\geq 4),\ 2^8\split\Orth^+_{8}(2),\ 2^6\split \Orth^-_{6}(2),\ 2^6\split \Sp_{6}(2).
\end{equation}
The parameters are listed in Table \ref{table4}.

{\small
\begin{table}[t]
\begin{center}
\renewcommand{\arraystretch}{1.5}
\begin{tabular}{|rc||c|c|c|c|c|c|}
\hline
$G$&&$\nu$&$k$&$s$&$g$&$c$&$d$\\
\hhline{|==#=|=|=|=|=|=|}
$\Sym_3$&&$3$&$2$&$-1$&$2$&$\frac{6}{5}$&$3$\\
\hline
$\Sym_n$&\hskip-.5em$n\geq 4$&$\frac{n(n-1)}{2}$&$2n-4$&$-2$&$\frac{n(n-3)}{2}$&$\frac{n(n-1)}{n+2}$&$\frac{n(n-1)}{2}$\\
\hline
$\EE\split\Sym_n$&\hskip-.5em$n\geq 4$&$n(n-1)$&$4n-8$&$-4$&$\frac{n(n-3)}{2}$&$n-1$&$n(n-1)$\\
\hline
$\EE^2\split\Sym_n$&\hskip-.5em$n\geq 4$&$2n(n-1)$&$8n-16$&$-8$&$\frac{n(n-3)}{2}$&$n$&$\frac{(3n-1)n}{2}$\\
\hhline{|==#=|=|=|=|=|=|}
$\Orth^-_6(2)$&&$36$&$20$&$-4$&$15$&$\frac{36}{7}$&$36$\\
\hline
$\Sp_6(2)$&&$63$&$32$&$-4$&$35$&$\frac{63}{10}$&$63$\\
\hline
$\Orth^+_8(2)$&&$120$&$56$&$-4$&$84$&$\frac{15}{2}$&$120$\\
\hline
$2^6\split\Orth^-_6(2)$&&$72$&$40$&$-8$&$15$&$6$&$57$\\
\hline
$2^6\split\Sp_6(2)$&&$126$&$64$&$-8$&$35$&$7$&$91$\\
\hline
$2^8\split\Orth^+_8(2)$&&$240$&$112$&$-8$&$84$&$8$&$156$\\
\hline
$\Orth^-_{8}(2)$&&$136$&$72$&$-8$&$51$&$\frac{34}{5}$&$85$\\
\hline
$\Sp_{8}(2)$&&$255$&$128$&$-8$&$135$&$\frac{15}{2}$&$120$\\
\hline
$\Orth^+_{10}(2)$&&$496$&$240$&$-8$&$340$&$8$&$156$\\
\hline
\end{tabular}
\end{center}
\caption{List of $(1/2,1/2)$-realizable groups}
\label{table4}
\end{table}
}
We have already seen in Proposition \ref{prop: root} that the groups in Table \ref{Table: root} are realizable by VOAs satisfying (I)--(IV).
On the other hand, the group $\Orth^+_{10}(2)$ is also realized by $V^+_{\sqrt{2}E_8}$ as shown by Griess \cite{Gr} (see Note \ref{note: Griess}) and the VOA realizing $\Sp_8(2)$ is constructed in \cite{KM}.
Since $\Orth^-_8(2)=\Orth^+_{10}(2)^{\Line}$, we can show by the same argument as in the proof of Proposition \ref{prop: root} using Lemma \ref{lemma: U2}
that this group is also realized by an appropriate subVOA of $V^+_{\sqrt{2}E_8}$.

\vskip4ex

Thus we have verified our main result:

\begin{theorem}
\sl
A centerfree 3-transposition group is $(1/2,1/2)$-realizable\ by a VOA satisfying {\rm (I)--(III)} if and only if it is the direct product of a finite number of groups of type listed in Table \ref{table4}.
\label{MainTheorem}
\end{theorem}
\section{Miscellaneous results and speculation}
\subsection{Groups $(1/2,1/2)$-realized by a VOA of class ${\cal S}^4$}
Let us consider the property of being of class ${\cal S}^4$ considered in \cite{Ma}.
Let $V$ be a VOA of class ${\cal S}^4$ satisfying (I)--(IV) with an idempotent $e\in B$ with the central charge $1/2$ such that $L(1/2,1/16)$-components are absent in the decomposition with respect to the Virasoro action generated by the vector $e$.
Then the dimension of the Griess algebra is determined by the rank $c=c^{}_V$  by\footnote{The formula (3.6) in \cite{Ma} should read  $(-22+2c)d=c(-37-10c)$.}
\begin{equation}
d=\frac{c(-37 - 10 c)}{-22 + 2 c}.
\end{equation}
Among such pairs of $c$ and $d$, the ones which fit $(1/2,1/2)$-realizable $3$-transposition groups are 
\begin{equation}
\Sym_2\subset \Sym_3\subset 2^4\split \Sym_4\subset \Orth^-_8(2)\subset \Sp_8(2)\subset \Orth^+_{10}(2),
\end{equation}
where we have included the case $\Sym_2$ additionally.
This sequence resembles the pattern of exceptional series $A_1\subset A_2\subset D_4\subset E_6\subset E_7\subset E_8$ of simply-laced root systems.
\subsection{$(1/2,1/16)$-realizable groups and the Monster}
Next, let us consider the case when $(G,D)$ is a centerfree 3-transposition group which is $(1/2,1/16)$-realizable.
Then, by Theorem \ref{least eigenvalue}, the least eigenvalue of the diagram is greater than or equal to $-64$.

Such groups are $(1/2,1/16)$-realizable at the level of the algebra $B$; The problem is to show the existence of a VOA satisfying (I)--(IV) whose Griess algebra coincides with $B$.

\vskip2ex

Now considering the moonshine module $V^\natural$, we have the following observation (cf.{} Note \ref{Miyamoto 3C}):
Suppose we have a subset $D$ of the Monster such that
\begin{enumerate}[\rm(a)]
\item
$D$ consists of involutions from the conjugacy class 2A.
\item
The product of any two elements of $D$ must fall in 1A, 2B or 3C.
\item
$D$ is stable under conjugation by elements of $D$.
\end{enumerate}
Then the pair $(G,D)$ where $G=\langle D\rangle$ must be a 3-transposition group for which the least eigenvalue of the diagram is greater than or equal to $-64$ and the central charge is less than or equal to $24$.

However, there is not much information on $(1/2,1/16)$-realizable groups; The author does not know the possibility of such a subgroup in the Monster except for $\Sym_3$.
\subsection{Simple group of order 168 and Steiner triple system}
Consider the projective plane over $\FF_2$, called the Fano plane, consisting 
of seven points $x_1,\ldots,x_7$ for which the lines are
$$\{x_1,x_2,x_3\}, \{x_1,x_4,x_5\}, \{x_1,x_6,x_7\}, \{x_2,x_4,x_6\}, 
\{x_2,x_5,x_7\}, \{x_3,x_4,x_7\}.$$
Though this is not a Fischer space, we may associate an algebra
 $\tildeb{B}$ with a bilinear form as in Table \ref{table1} in Subsect.\ref{Fischer space and Griess algebra} as 
well. Then the bilinear form is actually invariant with respect to 
the multiplication and the algebra has a unity if $\cw\ne -2/3$. 

However, this algebra does not satisfy the fusion property of 
binary type with respect to the idempotent corresponding to
 $\tildeb{x}_i$ for each $i=1,\ldots,7$, and the corresponding 
involution does not give us an automorphism of the algebra. 

Nevertheless, for $\ccc=4/5$ and $\cw=2/3$, the algebra $B$ associated with 
the Fano plane is isomorphic to the Griess algebra of a VOA of rank $14/5$ for 
which the automorphism group agrees with the automorphis group of the Fano 
plane: the simple group of order $168$ isomorphic to $\GL_3(2)$ and 
$\PSL_2(7)$. Such a VOA is in fact realized as a certain subVOA $V$ of 
$V^+_{\sqrt{2}D_4}$. Note that the VOA $V^+_{\sqrt{2}D_4}$ is isomorphic to 
the Hamming code VOA $V_{H_8}$, whose automorphism group has the shape 
$2^6\split (\GL_3(2)\times \Sym_3)$; see \cite{MM}. The VOA $V$ can 
be equally realized as a full subVOA of the VOA associated with the level one 
integrable highest weight representation of the affine Lie algebra of type 
$G^{(1)}_2$ by taking the fixed-point subspace with respect to the action of 
an elementary Abelian group $2^3$. 

\vskip2ex

Now, Let $X$ be any partial linear space $X$ such that each of the 
lines consists of three points.
Then we may associate the algebra $\tildeb{B}$ with the bilinear form by Table \ref{table1}.
Then obviously the multiplication is commutative and the form is symmetric.
Moreover, the form $(\ |\ )$ of the algebra $\tildeb{B}$ is 
invariant with respect to the multiplication $\cdot$ if and only 
if the space $X$ has the following properties:
\begin{enumerate}[\rm(1)]
\item
If $x\sim y\sim z$ then $\Col(x\compo y,z)=\Col(x,y\compo z)$.
\item
If $x\perp y\sim z\perp x$ then $x\perp y\compo z$.
\end{enumerate}
Here $\Col(x,y)=*$ if $x* y$ where $*$ is either $=$, $\perp$ or $\sim$.

For instance, any Steiner triple system clearly satisfies the condition (1) with the condition (2) being trivial.
It will be interesting to study the spaces $X$ satisfying the conditions (1) and (2) above.
%
%%%%%%%%%%%%%%%%%%%%%%%%%%%%%%%%%%%%%%%%%%%%%%%%%%%%%%%%%%%%%%%%%%%%%%%%%%%%%%
\bibliographystyle{unsrt}

\begin{thebibliography}{DMZ}
%
\bibitem[As]{As} M.~Aschbacher: 3-transposition groups. Cambridge tracts in Mathematics 124, Cambridge University Press, 1997. 
%
\bibitem[Bor]{Bor} R.E.~Borcherds: Vertex algebras, Kac-Moody algebras, and the monster, Proc.{} Nat'l.{} Acad.{} Sci.{} USA., {\bf 83}, (1986), 3068--3071.
%
\bibitem[Bos]{Bos} R.C.~Bose: Strongly regular graphs, partial geometries, and partially balanced designs, Pacific J.{} Math., {\bf 13}, (1963), 389-419.
%
\bibitem[CH]{CH} H.~Cuypers and J.I.~Hall: The $3$-transposition Groups with Trivial Center. J.{} Algebra, {\bf 178}, (1995), 149--193.
%
\bibitem[Co]{Co} J.H.~Conway: A simple construction for the Fischer-Griess monster group, Invent.{} Math., {\bf 79}, (1985), 513-540.
%
\bibitem[FHL]{FHL} I.~Frenkel, Y.-Z.~Huang and J.~Lepowsky: On axiomatic approaches to vertex operator algebras and modules. 
Memoirs Amer.{} Math.{} Soc.{} 104, 1993. American Mathematical Society.
%
\bibitem[DLMN]{DLMN} C.-Y.~Dong, H.-S.~Li, G.~Mason and S.P.~Norton: Associative subalgebras of the Griess algebra and related topics. The Monster and Lie algebras (Columbus, OH, 1996), 27--42, Ohio State Univ.{} Math.{} Res.{} Inst.{} Publ., 7, de Gruyter, Berlin, 1998. 
%
\bibitem[DMZ]{DMZ} C.-Y.~Dong, G.~Mason and Y.-C.~Zhu: Discrete series of the Virasoro algebra and the moonshine module. Algebraic groups and their generalizations: quantum and infinite-dimensional methods, 295--316, 
Proc.~Sympos.~Pure Math., 56, Part 2, 
Amer.~Math.~Soc., Providence, RI, 1994. 
%
\bibitem[Fi]{Fi} B.~Fischer: Finite Groups Generated by $3$-Transpositions. I. Invent.~Math., {\bf 13}, (1971), 232--246.
%
\bibitem[FLM]{FLM} I.B.~Frenkel, J.~Lepowsky and A.~Meurman: Vertex operator algebras and the Monster. Pure and Appl.~Math.~134, Academic Press, Boston, 1989.%
%
\bibitem[FQS]{FQS} D.~Friedan, Z.~Qiu and S.~Shenker: Details of non-unitarity proof for highest weight representations of the Virasoro algebra. Comm.{} Math.{} Phys.{} {\bf 107}, (1986), 535--542.
%
\bibitem[Gi]{Gi} P.~Ginsparg: Applied conformal field theory. Champs, cordes, et ph\'enom\`enes critiques, Les Houches Session XLIX, 1988, 3--168, Elsevier, 1989.
%
\bibitem[Gr]{Gr} R.L.~Griess, Jr: A vertex operator algebra related to $E\sb 8$ with automorphism group ${\rm O}\sp +(10,2)$. The Monster and Lie algebras, 43--58, Ohio State Univ.{} Math.{} Res.{} Inst.{} Publ., 7, de Gruyter, Berlin, 1998. 
%
\bibitem[Ha1]{Ha1} J.I.~Hall: Graphs, geometry, 3-transpositions, and symplectic $\FF_2$-transvection groups, Proc.~London Math.{} Soc.{} (Ser.~3), {\bf 58}, (1989), 89--111.
%
\bibitem[Ha2]{Ha2} J.I.~Hall: Some 3-transposition groups with non-central 2-subgroups, Proc.~London Math.~Soc.~(Ser.~3), {\bf 58}, (1989), 112--136.
%
\bibitem[Hi]{Hi} D.G.~Higman: Finite permutation groups of rank $3$, Math.~Z., {\bf 86}, (1964), 145--156.
%
\bibitem[KM]{KM} M.~Kitazume and M.~Miyamoto: 3-transposition automorphism groups of VOA. Groups and combinatorics -- in memory of Michio Suzuki, 315--324, Adv.{} Stud.{} Pure Math., 32, Mathematical Society of Japan, Tokyo, 2001. 
%
\bibitem[La]{La} C.-H.~Lam: Lattice vertex operator algebra $V_{\sqrt{2}E_8}$ and an algebra of Miyamoto of central charge $1/2+21/22$, in RIMS Kokyuroku 1327.
%
\bibitem[LY]{LY} C.-H.~Lam and H.~Yamada: $\ZZ_2\times\ZZ_2$ code VOA. J.~Algebra, {\bf 224}, (2000), 268--291.
%
\bibitem[MM]{MM} A.~Matsuo and M.~Matsuo: The automorphism group of the Hamming code vertex operator algebra, J.~Algebra, {\bf 228}, (2000), 204-226.
%
\bibitem[Ma]{Ma} A.~Matsuo: Norton's trace formulae for the Griess algebra of a vertex operator algebra with larger symmetry. Comm.~Math.~Phys.~{\bf 224}, no.~3, (2001), 565--591. 
%
\bibitem[MN]{MN} A.~Matsuo and K.~Nagatomo: Axioms for a vertex algebra and the locality of quantum fields, MSJ Memoirs 4. Mathematical Society of Japan, 1999.
%
\bibitem[Mi1]{Mi1} M.~Miyamoto: Griess algebras and conformal vectors in vertex operator algebras, J.~Algebra, {\bf 179}, (1996), 528--548.
%
\bibitem[Mi2]{Mi2} M.~Miyamoto: VOAs generated by two conformal vectors whose $\tau$-involutions generate $S_3$. Preprint.
%
\bibitem[No]{No} S.~Norton: On the group Fi$_{24}$. Geom.~Dedicata, {\bf 25}, (1988), 483--501.
%
\bibitem[Sh]{Sh} H.~Shimakura: Automorphism group of the vertex operator algebra $V^+_L$. Preprint.
%
\bibitem[We1]{We1} R.~Weiss: $3$-transpositions in infinite groups. Math.{} Proc.{} Cambridge Philos.{} Soc., {\bf 96}, no.~3, (1984), 371--377.
%
\bibitem[We2]{We2} R.~Weiss: A uniqueness lemma for groups generated by 3-transpositions. Math.{} Proc.{} Cambridge.{} Philos.{} Soc., {\bf 97}, (1985), 421--431.
%
\end{thebibliography}

\end{document}